\documentclass[12pt,draftcls,onecolumn]{IEEEtran}

\usepackage{amssymb}
\usepackage{epsfig}
\usepackage{psfrag}
\usepackage{latexsym}
\usepackage{amsmath}
\usepackage{amsfonts}
\usepackage{comment}
\usepackage{setspace}
\usepackage{epsfig}
\usepackage{dsfont}
\usepackage{txfonts}
\title{\LARGE \bf
$\mc{H}_2$-Optimal Decentralized Control over Posets: \\A State-Space Solution for State-Feedback}

\author{Parikshit Shah and Pablo A. Parrilo  \\
MIT}

\def\T={\buildrel {\scriptscriptstyle\triangle} \over =}

\def\mc{\mathcal}
\def\mb{\mathbb}

\def\ba{\begin{array}}
\def\ea{\end{array}}
\def\ll{\left[}
\def\rr{\right]}
\def\l{\left\{}
\def\r{\right\} }
\def\down{{\downarrow}}
\def\up{{\uparrow}}
\def\inc{\mathcal{I}(\mathcal{P})}
\newtheorem{theorem}{Theorem}

\newtheorem{cor}{Corollary}
\newtheorem{definition}{\indent Definition}

\newtheorem{lemma}{Lemma}
\newtheorem{example}{Example}
\newtheorem{assump}{Assumption}
\newenvironment{remark}[1][Remark]{\begin{trivlist}
\item[\hskip \labelsep {\bfseries #1}]}{\end{trivlist}}
\newenvironment{remarks}[1][Remarks]{\begin{trivlist}
\item[\hskip \labelsep {\bfseries #1}]}{\end{trivlist}}

\begin{document}

\maketitle
\thispagestyle{plain}
\pagestyle{plain}

\bibliographystyle{plain}

\begin{abstract}
We develop a complete state-space solution to $\mc{H}_2$-optimal decentralized
control of poset-causal systems with state-feedback. Our solution is based
on the exploitation of a key separability property of the problem, that
enables an efficient computation of the optimal controller by solving a
small number of uncoupled standard Riccati equations. Our approach gives
important insight into the structure of optimal controllers, such as controller degree
bounds that depend on the structure of the poset. A novel element in our
state-space characterization of the controller is a remarkable pair of
transfer functions, that belong to the incidence algebra of the poset, are
inverses of each other, and are intimately related to prediction of the
state along the different paths on the poset. The results are illustrated by
a numerical example.

\end{abstract}

\section{Introduction}\label{sec:1}
Finding computationally efficient algorithms to design decentralized controllers is a challenging area of research (see e.g. \cite{poset,blondel} and the references therein). Current research suggests that while the problem is hard in general, certain classes with special information structures are tractable via convex optimization techniques. In past work, the authors have argued that communication structures modeled by partially ordered sets (or \emph{posets}) provide a rich class of decentralized control systems (which we call \emph{poset-causal} systems) that are amenable to such an approach \cite{poset}. Posets have appeared in the control theory literature earlier in the context of team theory \cite{Ho}, and specific posets (chains) have been studied in the context of decentralized control \cite{petros2}. Poset-causal systems are also related to the class of systems studied more classically in the context of \emph{hierarchical systems} \cite{Mesarovic}, \cite{Findeisen}, where abstract notions of hierarchical organization of large-scale systems were introduced and their merits were argued for.

While it is possible to design optimal decentralized controllers for a fairly large class of systems known as quadratically invariant systems in the frequency domain via the Youla parametrization \cite{quadinv}, there are some important drawbacks with such an approach. Typically Youla domain techniques are not computationally efficient, and the degree of optimal controllers synthesized with such techniques is not always well-behaved. In addition to computational efficiency, issues related to numerical stability also arise. Typically, operations at the transfer function level are inherently less stable numerically. Moreover, such approaches typically do not provide insight into the structure of the optimal controller. These drawbacks emphasize the need for \emph{state-space techniques} to synthesize optimal decentralized controllers. State-space techniques are usually computationally efficient, numerically stable, and provide degree bounds for optimal controllers. In our case we will also show that the solution provides important insight into the structure of the controller.

In this paper we consider the problem of designing $\mc{H}_2$ optimal decentralized controllers for poset-causal systems. The control objective is the design of optimal feedback laws that have access to local state information. We emphasize here that different subsystems \emph{do not} have access to the global state, but only the local states of the systems in a sense that will be made precise in the next section. The main contributions in the paper are as follows:
\begin{itemize}
\item We show a certain crucial separability property of the problem under consideration. This result is outlined in Theorem~\ref{theorem:2}. This makes it possible to decompose the decentralized control problem over posets into a collection of standard centralized control problems.
\item We give an explicit state-space solution procedure in Theorem~\ref{theorem:3}. To construct the solution, one needs to solve standard Riccati equations (corresponding to the different sub-problems). Using the solutions of these Riccati equations, one constructs certain block matrices and provides a state-space realization of the controller.
\item We provide bounds on the degree of the optimal controller in terms of a parameter $\sigma_{\mc{P}}$ that depends only on the order-theoretic structure of the poset (Corollary \ref{cor:degbounds}).
\item In Theorem~\ref{theorem:4} we briefly describe the structural form of the optimal controller. We introduce a novel pair of transfer functions $(\Phi, \Gamma)$ which are inverses of each other, and which capture the prediction structure in the optimal controller. We call $\Phi$ the \emph{propagation filter}, it plays a role in propagating local signals (such as states) upstream based on local information. We call $\Gamma$  the \emph{differential filter}, it corresponds to computation of differential improvement in the prediction of the state at different subsystems. The discussion related to structural aspects is brief and informal in this paper and have been formalized in the paper \cite{poset_struct}. 
\item We state a new and intuitive decomposition of the structure of the optimal controller into certain local control laws.
\end{itemize}
\subsection{Related Work}
It is well-known that in general decentralized control is a hard
problem, and significant research efforts have been directed towards
its many different aspects; see for instance the classical survey
\cite{sandell} for many of the earlier results. More recently, Blondel
and Tsitsiklis \cite{Blondel2} have shown that in certain instances,
decentralized control problems are computationally
intractable, in particular they show that the problem of finding bounded-norm, block-diagonal stabilizing controllers in the presence of output-feedback is NP-hard. On the other hand, Voulgaris \cite{petros2}, \cite{petros1}
presented several cases where decentralized control problems are amenable to a convex reparametrization and therefore computationally tractable. Lall and Rotkowitz generalize these ideas in terms of a property called quadratic invariance \cite{quadinv}, we discuss connections to their work later. In past work \cite{poset}, we have shown that posets provide a unifying umbrella to describe these tractable examples under an appealing
theoretical framework. 

Partially ordered sets (posets) are very well studied objects in combinatorics. The associated notions of incidence algebras and Galois connections were first studied
by Rota \cite{rota} in a combinatorics setting. Since then,
order-theoretic concepts have been used in engineering and computer
science; we mention a few specific works below. In control theory,
ideas from order theory have been used in different ways. Ho and Chu used posets to study team
theory problems \cite{Ho}. They were interested in sequential decision
making problems where agents must make decisions at different time
steps. They study computational and structural properties of optimal decision policies when the
problems have poset structure. Mullans and Elliot \cite{mullans} use
posets to model the notions of time and causality, and study
evolution of systems on locally finite posets. Wyman \cite{Wyman} has studied time-varying linear-systems evolving on locally finite posets in an algebraic framework, including aspects related to realization theory and duality.  In computer science, Cousot and Cousot used
these ideas to develop tools for formal verification of computer
programs in their seminal paper \cite{cousot}.   Del Vecchio and Murray
\cite{delvecchio} have used ideas from lattice and order theory to
construct estimators for continuous states in hybrid systems.

More recently, the authors of this paper have initiated a systematic study of decentralized control problems from the point of view of partial order theory. In \cite{poset}, we introduce the partial order framework and show how several well-known classes of problems such as nested systems \cite{petros2} fit into the partial order framework. In \cite{spatinv}, we extend this poset framework to spatio-temporal systems and generalize certain results related to the so called ``funnel causal systems'' of Bamieh et. al \cite{Bamieh}. In \cite{poset2}, we show that a class of time-delayed systems known to be amenable to convex reparametrization \cite{quadinv} also has an underlying poset structure. In that paper, we also study the close connections between posets and another class of decentralized control problems known as quadratically invariant problems. While this poset framework provides a lens to view all these examples in a common intuitive framework, a systematic study of state-space approaches has been lacking.

In an interesting paper by Swigart and Lall \cite{Swigart}, the authors consider a state-space approach to the $\mc{H}_2$ optimal controller synthesis problem over a particular poset with two nodes. Their approach is restricted to the finite time horizon setting (although in a subsequent paper \cite{Swigart2}, they extend this to the infinite time horizon setting), and uses a particular decomposition of certain optimality conditions. In this setting, they synthesize optimal controllers and provide insight into the structure of the optimal controller. These results are also summarized in the thesis \cite{swigart_thesis}. By using our new separability condition (which is related to their decomposition property, but which we believe to be more fundamental) we significantly generalize those results in this paper. We provide a solution for \emph{all posets} and for the \emph{infinite time horizon}. In recent work \cite{quadinv}, Rotkowitz and Lall proposed a state-space technique to solve $\mc{H}_2$ optimal control problems for quadratically invariant systems (which could be used for poset-causal systems). 
An important drawback of their reformulation is that one would need to solve larger Riccati equations. Our approach for poset-causal systems is more efficient computationally. Moreover, our approach also provides insight into the structure of the optimal controllers.
 
The rest of this paper is organized as follows. In Section~\ref{sec:2} we introduce the necessary preliminaries regarding posets, the control theoretic framework and notation. In Section~\ref{sec:3} we describe our solution strategy. In Section~\ref{sec:4} we present the main results. We devote Section~\ref{sec:5} to a discussion of the main results, and their illustration via examples. 

\section{Preliminaries}\label{sec:2}
In this section we introduce some concepts from order theory. Most of these concepts are well studied and fairly standard, we refer the reader to \cite{aigner},\cite{davey} for details.
\subsection{Posets}
\begin{definition}
A partially ordered set (or \emph{poset}) $\mc{P}=(P, \preceq)$ consists of a set $P$
along with a binary relation $\preceq$ which has the following properties:
\begin{enumerate}
 \item $a \preceq a$ (reflexivity),
\item  $a \preceq b$ and $b \preceq a$ implies $a=b$ (antisymmetry), 
\item $a \preceq b$ and $b \preceq c$ implies $a \preceq c$ (transitivity).
\end{enumerate}
\end{definition}
We will sometimes use the notation $a \prec b$ to denote the strict order relation $a \preceq b$ but $a \neq b$.

In this paper we will deal with finite posets (i.e. $|P|$ is finite).
It is possible to represent a poset graphically via a \emph{Hasse diagram} by representing the transitive reduction of the poset as a graph (i.e. by drawing only the minimal order relations graphically, a downward arrow representing  the relation $\preceq$, with the remaining order relations being implied by transitivity).
\begin{example} \label{example:1}
An example of a poset with three elements (i.e., $P=\left\{ 1,
2, 3 \right\}$) with order relations $1 \preceq 2$ and $1 \preceq 3 $
is shown in Figure \ref{fig:2.1}(b).
\begin{figure}[htbp]
  \begin{center}
  \includegraphics[scale=0.6]{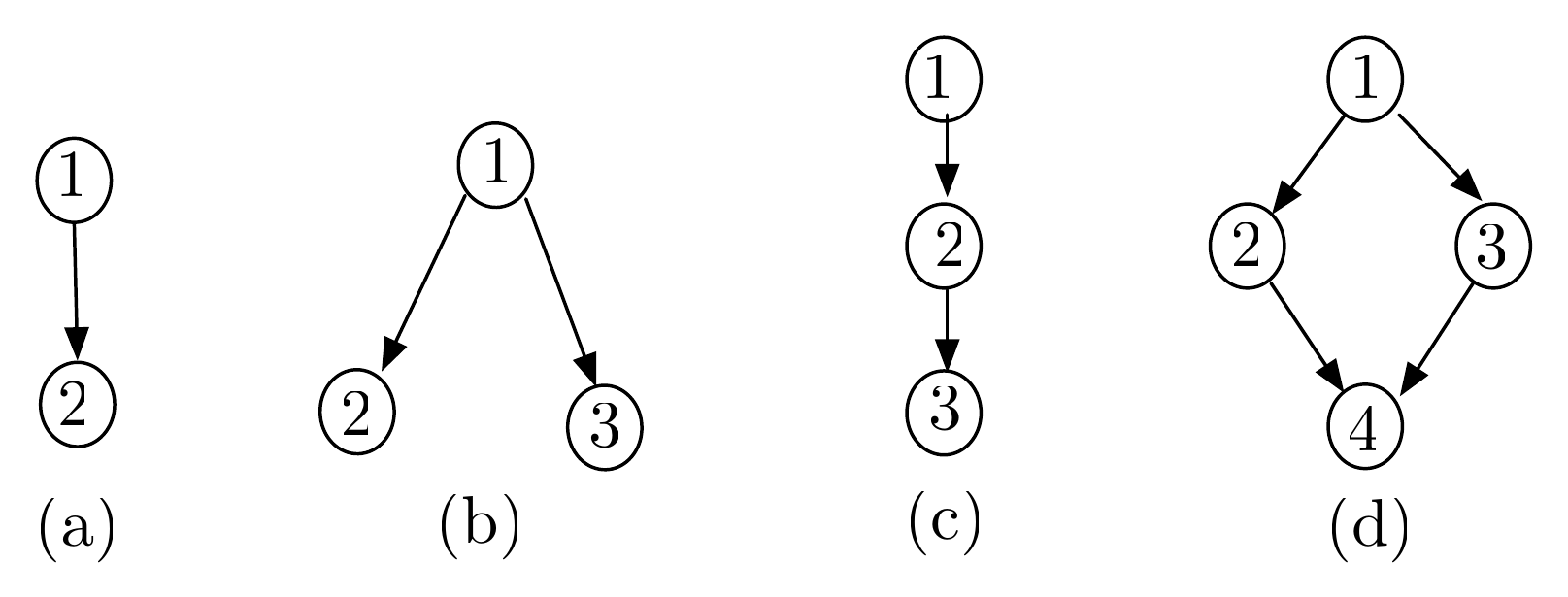}
  \end{center}
  \caption{Hasse diagrams of some posets.}
  \label{fig:2.1}
\end{figure}
\end{example}
Let $\mc{P}=(P, \preceq)$ be a poset and let $p\in P$. We define $\down p=\l q\in P \; |\; p \preceq q \r$ (we call this the \emph{downstream set}).
\footnote[3]{We have reversed conventions with respect to some of our conference papers, wherein the Hasse diagrams are drawn with upward arrows and the set $\up p$ corresponds to the set $\l q\in P \; |\; p \preceq q \r$. The present convention has been adopted to make the presentation more intuitive. For example the downstream set at $p$ corresponds to elements drawn lower in the Hasse diagram. It also corresponds to the elements that are ``in the future'' with respect to $p$, in keeping with the intuition that information in a river propagates ``downstream''.}
Let $\down \down p = \l q\in P \; | \; p \preceq q, q \neq p \r$. 
Similarly, let  $\up p=\l q\in P \; | \; q \preceq p \r$ (called the \emph{upstream set}), and $\up \up p= \l q\in P \; | \; q \preceq p, q \neq p \r$. We define $ \down \up p = \l q\in P \; | \; q \npreceq p, q \npreceq p \r$ (called the \emph{off-stream set}); this is the set of \emph {uncomparable} elements that have no order relation with respect to $p$. 
Define an \emph{interval} $[i,j] =\l p \in P \; | \; i \preceq p \preceq j \r$.
 A \emph{minimal element} of the poset is an element $p \in P$ such that if $q \preceq p$ for some $q \in P$ then $q=p$. (A maximal element is defined analogously).
 
 In the poset shown in Figure \ref{fig:2.1}(d), $\down 1 = \l 1, 2, 3, 4 \r$, whereas $\down \down 1 = \l  2, 3, 4 \r$. Similarly $\up \up 1 = \emptyset$, $\up 4 = \l 1, 2, 3, 4 \r$, and $\up \up 4 = \l 1, 2, 3 \r$. The set $\down \up 2=\l 3 \r$. 

\begin{definition} \label{def:inc_algebra}
Let $\mc{P}=(P,\preceq)$ be a poset. Let $\mathbb{Q}$ be a ring. The set of all functions
$f:P \times P \rightarrow \mathbb{Q}$
with the property that $f(x,y)=0$ if $y \npreceq x$ is called the
\emph{incidence algebra} of $\mc{P}$ over $\mathbb{Q}$. It is denoted by
$\mc{I}(\mc{P})$.
\footnote[1]{Standard definitions of the incidence algebra use an opposite convention, namely $f(x,y)=0$ if $x \npreceq y$. Thus, the matrix representation of the incidence algebra is typically a transposal of the matrix representations that appear here. For example, while the incidence algebra of a chain is the set of lower-triangular matrices in this paper, in standard treatments it would appear as upper-triangular matrices. We reverse the convention so that in a control theoretic setting one may interpret such matrices as representing \emph{poset-causal} maps.}

\end{definition}
When the poset $\mc{P}$ is finite, the elements in the incidence
algebra may be thought of as matrices with a specific sparsity
pattern given by the order relations of the poset in the following way. 
One indexes the rows and columns of the matrices by the elements of $P$. Then the $(i,j)$ entry for $i,j \in P$ of the matrix corresponds to $f(i,j)$. By Definition \ref{def:inc_algebra}, if $j \npreceq i$ then the $(i,j)$ entry of the matrix must be zero.
An example of an element of $\mc{I}(\mc{P})$ for the poset from
Example 1 (Fig. \ref{fig:2.1}(b)) is:
\begin{equation*}
\zeta_{\mc{P}}=\left[ \begin{array}{ccc} 1 & 0 & 0 \\
1 & 1 & 0 \\
1 & 0 & 1 \end{array} \right].
\end{equation*}
Given two functions $f, g \in \mc{I}(\mc{P})$, their sum
$f+g$ and scalar multiplication $c f$ are defined as usual. The
product $h=f \cdot g$ is defined by
$h(x,y)=\sum_{z \in P}f(x,z)g(z,y)$.
Note that the
above definition of function multiplication is made so that it is
consistent with standard matrix multiplication.
\begin{lemma} \label{lemma:inc}
Let $\mc{P}$ be a poset. Under the usual definition of addition and
multiplication as defined in (1) the incidence algebra is an associative algebra
(i.e. it is closed under addition, scalar multiplication and
function multiplication).
\end{lemma}
\begin{proof}
The proof is standard, see for example \cite{poset}.
\end{proof}
 Given $i \preceq j$, let $[i \rightarrow j]$ denote the set of all chains from $i$ to $j$ of the form $\{i, i_1\}, \ldots, \{i_k, j\}$ such that $i \preceq i_1 \preceq \cdots \preceq i_k\preceq j$. 
 For example, in the poset in Fig. \ref{fig:2.1}(c), $[1 \rightarrow 3]=\l \l \l1,2\r, \l 2,3\r \r, \l1,3 \r \r$. 
A standard corollary of Lemma \ref{lemma:inc} is the following.
\begin{cor} \label{cor:inc}
Suppose $A \in \mc{I}(\mc{P})$. Then $A$ is invertible if and only if $A_{ii}$ is invertible for all $i \in P$. Furthermore $A^{-1} \in \mc{I}(\mc{P})$, and the inverse is given by:
$$
[A^{-1}]_{ij}= \left\{ \ba{ll} A_{ii}^{-1}\sum_{p_{ij}\in[j \rightarrow i]} \prod_{\{l,k \} \in p_{ij}}(-A_{lk}A_{kk}^{-1}) &\text{ if } i \neq j\\
A_{ii}^{-1} & \text{ if }  i= j. \ea \right.
$$
\end{cor}
\subsection{Control Theoretic Preliminaries}
We consider the following state-space system in continuous time:

\begin{equation} \label{eq:2}
\begin{split}
\dot{x}(t)&=Ax(t)+Fw(t)+Bu(t) \\
z(t)&=Cx(t)+Du(t).
\end{split}
\end{equation}
In this paper we present the continuous time case only, however, we wish to emphasize that analogous results hold in discrete time in a straightforward manner.
In this paper we consider what we will call \emph{poset-causal systems}. We think of the system matrices $(A,B,C,D,F)$ to be partitioned into blocks in the following natural way. Let $\mc{P}=(P,\preceq)$ be a poset with $P=\left\{1,\ldots, p \right\}$. 
We think of this system as being divided into $p$ sub-systems, with sub-system $i$ having some states $x_i(t) \in \mb{R}^{n_i}$, and control inputs $u_i(t) \in \mb{R}^{m_i}$ for $i\in \left\{1,\ldots, p \right\}$.  The external output is $z(t) \in \mb{R}^{l}$. The signal $w(t)$ is a disturbance signal. (To use certain standard state-space factorization results, we assume that $C^{T}D=0$ and $D^{T}D \succ 0$, these assumptions can be relaxed in a straightforward way). The states and inputs are partitioned in the natural way such that the sub-systems correspond to elements of the poset $\mc{P}$ with $x(t)=\left[ x_1(t)\left| x_2(t) \left| \ldots \left| x_p(t) \right. \right. \right. \right] ^{T}$, and $u(t)=\left[ u_1(t)\left| u_2(t) \left| \ldots \left| u_p(t) \right. \right. \right. \right] ^{T}$. This naturally partitions the matrices $A, B, C, D, F$ into appropriate blocks so that $A=\left[A_{ij} \right]_{i,j \in P}$,  $B=\left[B_{ij} \right]_{i,j \in P}$,  $C=\left[C_{j} \right]_{j \in P}$ (partitioned into columns),  $D=\left[D_{j} \right]_{j \in P}$, $F=\left[ F_{ij} \right]_{i,j \in P}$. (We will throughout deal with matrices at this block-matrix level, so that $A_{ij}$ will unambiguously mean the $(i,j)$ block of the matrix $A$.)
Using these block partitions, one can define the incidence algebra at the block matrix level in the natural way.
We denote by $\mc{I}_{A}(\mc{P}), \mc{I}_B(\mc{P})$ the block incidence algebras corresponding to the block partitions of $A$ and $B$. 

We will further assume that $F$ is block diagonal and full column rank.
\footnote[2]{More generally we can assume that for the system under consideration \eqref{eq:2}, $F \in \inc$ and the diagonal blocks are full column rank. Operating under this assumption, one can perform an invertible coordinate transformation $T \in \inc$ on the states so that $T^{-1}F$ is block diagonal. Since $T$ can be chosen to be in the incidence algebra, $T^{-1}AT,T^{-1} B \in \inc$. Hence, without loss of generality, we assume that $F$ is block diagonal.}
Often, matrices will have different (but compatible) dimensions and the block structure will be clear from the context. In these cases, we will abuse notation and will drop the subscript and simply write $\mc{I}(\mc{P})$.

\begin{definition}
We say that a state-space system is $\mc{P}$-\emph{poset-causal} (or simply poset-causal) if $A \in \mc{I}_A(\mc{P})$ and   $B \in \mc{I}_B(\mc{P})$.
\end{definition}
\begin{example}\label{example:2}
We use this example to illustrate ideas and concepts throughout this paper. Consider the system
\begin{align*}
\dot{x}(t)&=Ax(t)+Fw(t)+Bu(t)\\
z(t)&=Cx(t)+Du(t) \\
y(t)&=x(t) ,
\end{align*}
with matrices
\begin{align} \label{eq:comp_example}
A&=\ll \ba{cccc}
-0.5     &    0     &    0    &     0 \\
   -1&  -0.25      &   0    &     0 \\
   -1   &      0  & -0.2   &      0 \\
   -1   & -1&  -1  & -0.1
 \ea \rr &
 B&= \ll \ba{cccc}
    1  &      0   &      0   &      0 \\
    1  &  1    &     0   &      0 \\
   1   &      0  &  1    &     0 \\
    1  &  1  & 1  & 1 \ea \rr \\
    C&=\ll \ba{c} I^{4 \times 4} \\
    		      0^{4 \times 4} \ea \rr \qquad \qquad \qquad F=I &
    D&=\ll \ba{c} 0^{4 \times 4} \\
    		      I^{4 \times 4} \ea \rr.
\end{align}
This system is poset-causal with the underlying poset described in Fig. \ref{fig:2.1}(d).
Note that in this system, each subsystem has a single input, a single output and a single state. The matrices $A$ and $B$ are in the incidence algebra of the poset. Furthermore, $F=I$.
\end{example}

Recall that the standard notion of causality in systems theory is based crucially on an underlying totally ordered index set (time). Systems (in LTI theory these are described by impulse responses) are said to be \emph{causal} if the support of the impulse response is consistent with the ordering of the index set: an impulse at time zero is only allowed to propagate in the increasing direction with respect to the ordering. This notion of causality can be readily extended to situations where the underlying index set is only \emph{partially ordered}. Indeed this abstract setup has been studied by Mullans and Elliott \cite{mullans}, and an interesting algebraic theory of systems has been developed. 

Our notion of poset-causality is very much in the same spirit. We call such systems poset-causal due to the following analogous property among the sub-systems. If an input is applied to sub-system $i$ via $u_i$ at some time $t$, the effect of the input is seen by the states $x_j$ for all sub-systems $j \in \down i$ (at or after time $t$). Thus $ \down i$ may be seen as the cone of influence of input $i$. We refer to this causality-like property as \emph{poset-causality}. This notion of causality enforces (in addition to causality with respect to time), causality with respect to the subsystems via a poset.
For most of this paper we will deal with systems that are poset-causal (with respect to some arbitrary but fixed finite poset $\mc{P}$). Before we turn to the problem of optimal control we state an important result regarding stabilizability of poset-causal systems by poset-causal controllers. \begin{theorem} \label{theorem:stab}
The poset-causal system \eqref{eq:2} is stabilizable by a poset-causal controller $K\in\mc{I}(\mc{P})$ if and only if the $(A_{ii},B_{ii})$ are stabilizable  for all $i \in P$. 
\end{theorem}
\begin{proof}
See Appendix.
\end{proof}
In this paper, we make the following important assumption about the stabilizability of the sub-systems. By the preceding theorem, this assumption is necessary and sufficient to ensure that the systems under consideration have feasible controllers. 
\begin{assump}
Given the poset-causal system of the form \eqref{eq:2}, we assume that the sub-systems $(A_{ii}, B_{ii})$ are stabilizable for all $i \in \left\{1, \ldots, p\right\}$. 
\end{assump}
In the absence of this assumption, there is no poset-causal stabilizing controller, and hence the problem of finding an optimal one becomes vacuous. This assumption is thus necessary and sufficient for the problem to be well-posed.
Moreover, in what follows, we will need the solution of certain standard Riccati equations. Assumption~1 ensures that all of these Riccati equations have well-defined stabilizing solutions. This stabilizing property of the Riccati solutions will be useful for proving internal stability of the closed loop system.

The system \eqref{eq:2} may be viewed as a map from the inputs $w, u$ to outputs $z,x$ via
\begin{align*}
z&=P_{11}w+P_{12}u \\
x&=P_{21}w+P_{22}u
\end{align*}
where \small
\begin{equation} \label{eq:1}
\ba{rl}
\ll
\ba{cc} P_{11} & P_{12} \\
P_{21} & P_{22}
\ea \rr&=\ll \ba{cc}
C(sI-A)^{-1}F & C(sI-A)^{-1}B+D \\
(sI-A)^{-1}F & (sI-A)^{-1}B
\ea\rr \\
&=\ll \ba{c|cc} 
A&F&B \\ \hline C & 0 & D \\ I & 0 & 0
\ea \rr.
\ea
\end{equation}
\normalsize
A controller $u=Kx$ induces a map $T_{zw}$ from the disturbance input $w$ to the exogenous output $z$ via 
$$
T_{zw}=P_{11}+P_{12}K(I-P_{22}K)^{-1}P_{21}.
$$
Thus, after the controller is interconnected with the system, the closed-loop map is $T_{zw}$. The objective function of interest is to minimize the $\mc{H}_2$ norm  \cite{Zhou} of $T_{zw}$ which we denote by $\| T_{zw}\|$.

\subsection{Information Constraints on the Controller} 
Given the system \eqref{eq:2}, we are interested in designing a controller $K$ that meets certain specifications. In traditional control problems, one requires $K$ to be proper, causal and stabilizing. One can impose additional constraints on the controller, for example require it to belong to some subspace. Such seemingly mild requirements can actually make the problem significantly more challenging.
This paper focuses on addressing the challenge posed by subspace constraints arising from particular decentralization structures. The decentralization constraint of interest in this paper is one where the controller mirrors the structure of the plant, and is therefore also in the block incidence algebra $\mc{I}_{K}(\mc{P})$ (we will henceforth drop the subscripts and simply refer to the incidence algebra $\mc{I}(\mc{P})$). This translates into the requirement that input $u_i$ only has access to $x_j$ for $j \in \up i$ thereby enforcing poset-causality constraints also on the controller. In this sense the controller has access to \emph{local} states, and we thus refer to it as a \emph{decentralized} state-feedback controller.
\subsection{Problem Statement}
Given the poset-causal system \eqref{eq:1} with poset $\mc{P}=(P, \preceq)$, $|P|=p$, solve the optimization problem:
\begin{equation} \label{eq:chap4_opt_controll}
\ba{rl}
\underset{K}{\text{minimize}} \text{ \ \ }  & \|P_{11}+P_{12}K(I-P_{22}K)^{-1}P_{21} \|^{2} \\
\text{subject to \ \ } & K\in \mc{I}(\mc{P}) \\
& K \text{ stabilizing}.
\ea
\end{equation} 
The main problem under consideration is to solve the above stated optimal control problem in the controller variable $K$. The feasible set is the set of all rational proper transfer function matrices that internally stabilize the system \eqref{eq:2}. In the absence of the decentralization constraints $K \in \mc{I}(\mc{P})$ this is a standard, well-studied control problem that has an efficient finite-dimensional \emph{state-space} solution \cite{Zhou}. The main objective of this paper is to construct such a solution for the poset-causal case.

\subsection{Notation}
Given a matrix $Q$, let $Q(j)$ denote the $j^{\text{th}}$ column of $Q$. We denote the $i^{th}$ component of the vector $Q(j)$ to be $Q(j)_i$. For a poset $\mc{P}$ with incidence algebra $\mc{I}(\mc{P})$, if $M \in \inc$ then recall that $M$ is sparse, i.e. has a zero pattern given by $M_{ij}=0$ if $j \npreceq i$. We denote the sparsity pattern of the $j^{th}$ column of the matrices in  $\mc{I}(\mc{P})$ by $\mc{I}(\mc{P})^{j}$. More precisely if $v$ is a vector (of length equal to that of the matrix $M$) whose components are denoted by $v_i$ then
$$
 \mc{I}(\mc{P})^{j}:=\left\{ v| v_i=0 \text{ for } j \npreceq i \right\}.
$$
In the above definition $v$ is understood to be a vector composed of $|P|$ blocks, with sparsity being enforced at the block level.

Given the data $(A,B,C,D)$, we will often need to consider sub-matrices or embed a sub-matrix into a full dimensional matrix by zero padding. Some notation for that purpose we will use is the following: 
\begin{enumerate}
\item Define $Q^{\down j}=[Q_{ij}]_{i \in \down j}$ (so that it is the $j^\text{th}$ column shortened to include only the nonzero entries).  
\item Also define $A( \down j)=[A(i)]_{i \in \down j}$ so that it is the sub-matrix of $A$ containing exactly those columns corresponding to the set $\down j$. 
\item Define $A( \down j, \down j)=[A_{kl}]_{k, l \in \down j}$ so that it is the $(\down j, \down j)$ sub-matrix of $A$ (containing exactly those rows and columns corresponding to the set $\down j$). 
\item Sometimes, given a block $|\down j|  \times |\down j|$ matrix we will need to embed it into a block matrix indexed by the original poset (i.e. a $p \times p$ matrix) by padding it with zeroes. Given $K$ (a block $|\down j|  \times |\down j|$ matrix) we define:
 $$
 [\hat{K}]_{l.m}=\left\{ \ba{l} K_{lm} \text{ if $l, m \in \down j$} \\ 0 \text{ otherwise. }
 \ea \right.
 $$  
\item $E_i=[\ba{ccccc} 0 & \ldots& I & \ldots& 0 \ea]^{T}$ be the tall block matrix (indexed with the elements of the poset) with an identity in the $i^{th}$ block row. 
\item Let $S \subseteq P$. Define $E_S=[E_i]_{i \in S}$. Note that given a block $p \times p$ matrix $M$, $ME_{\down j}=M(\down j)$ is a matrix containing the columns indexed by $\down j$. 
\item Given matrices $A_i, i \in P$, we define the block diagonal matrix:
$$
\text{diag}(A_i)=\ll \ba{ccc} 
A_{1} &  &  \\
& \ddots & \\
& &A_{p}
\ea \rr.
$$
\end{enumerate}
Recall that every poset $\mc{P}$ has a linear extension (i.e. a total order on $P$ that is consistent with the partial order $\preceq$). For convenience, we fix such a linear extension of $\mc{P}$, and all indexing of our matrices throughout the paper will be consistent with this linear extension (so that elements of the incidence algebra are lower triangular).

\begin{example}
Let $\mc{P}$ be the poset shown in Fig. \ref{fig:2.1}(d). We continue with Example \ref{example:2} to illustrate notation. (Note that $\down 2 = \l 2, 4 \r$). As per the notation defined above,
\begin{align*}
A^{\down 2}&=\ll \ba{c} -.25 \\ -1 \ea \rr & 
A(\down 2)&=\ll \ba{cc} 0 & 0 \\ -0.25 & 0 \\ 0 & 0 \\ -1 & -0.1 \ea \rr &
A(\down 2, \down 2)&=\ll \ba{cc}  -0.25 & 0 \\ -1 & -0.1 \ea \rr. \end{align*}
Also, if 
$K(\down 2, \down 2)=\ll \ba{cc} 1 & 2 \\ 3 & 4 \ea \rr, \text{then \ } 
\hat{K}(\down 2, \down 2)=\ll \ba{cccc} 0 & 0 & 0 & 0\\ 0 &1 & 0 & 2 \\0 & 0 & 0 & 0\\0 & 3 &0 & 4 \ea \rr$.

\end{example}
\section{Solution Strategy}\label{sec:3}
In this section we first remind the reader of a standard reparametrization of the problem known as the Youla parametrization. Using this reparametrization, we illustrate the main technical idea of this paper using an example.
\subsection{Reparametrization}
Problem \eqref{eq:chap4_opt_controll} as stated has a nonconvex objective function. Typically \cite{quadinv,poset}, this is convexified by a bijective change of parameters given by $R:=K(I-P_{22}K)^{-1}$ (though one typically needs to make a stability or prestabilization assumption). 
When the sparsity constraints are poset-causal (or quadratically invariant, more generally), this change of parameters preserves the sparsity constraints, and $R$ inherits the sparsity constraints of $K$. The resulting infinite-dimensional problem is convex in $R$. 

For poset-causal systems with state-feedback we will use a slightly different parametrization. 
We note that for poset-causal systems, the matrices $A$ and $B$ are both in the block incidence algebra. As a consequence of \eqref{eq:1},  ${P}_{21}$ and $P_{22}$ are also in the incidence algebra. This structure, which follows from the closure properties of an incidence algebra, will be extensively used. Since ${P}_{21}, P_{22} \in \mc{I}(\mc{P})$ the optimization problem \eqref{eq:chap4_opt_controll} maybe be reparametrized as follows. Set 
\begin{equation} \label{eq:param}
Q:=K(I-P_{22}K)^{-1}{P}_{21}.
\end{equation}
 Note that $P_{21}$ is left invertible, and a left inverse is given by
 $$
 P_{21}^{\dagger}=F^{\dagger}(sI-A),
 $$
 where $F^\dagger$ is the pseudoinverse of $F$, so that $F^\dagger F=I$ (note also that the pseudoinverse is block diagonal and hence in $\inc$).
 As a consequence, 
given $Q$, $K$ can be recovered using 
\begin{equation} \label{eq:5} 
K=Q{P}_{21}^{\dagger}(I+P_{22}Q{P}_{21}^{\dagger})^{-1}.
\end{equation} 
Since $I, P_{21}, {P}_{21}^{\dagger}, P_{22}$ all lie in the incidence algebra, $K \in \mc{I}(\mc{P})$ if and only if $Q \in \mc{I}(\mc{P})$.  
Using this reparametrization the optimization problem \eqref{eq:chap4_opt_controll} can be relaxed to:
\begin{equation}\label{eq:6a}
\ba{rl}
\underset{Q}{\text{minimize}} & \|P_{11}+P_{12}Q \|^{2} \\
\text{subject to} & Q \in \mc{I}(\mc{P}). \\
\ea
\end{equation}
\begin{remarks}
\begin{enumerate}
\item We note that ${P}_{21}^{\dagger}$, and hence \eqref{eq:5} may potentially be improper. However, we will \emph{prove} that for the optimal $Q$ in \eqref{eq:6a}, this expression is proper and corresponds to a rational controller $K^{*} \in \inc$. 
\item For the objective function to be bounded, the optimal $Q$ would have to render $P_{11}+P_{12}Q$ stable. However, one also requires that the overall system is internally stable. We relax this requirement on $Q$ and later show that $K^{*}$ is nevertheless internally stabilizing. Thus \eqref{eq:6a} is in fact a \emph{relaxation} of \eqref{eq:chap4_opt_controll}. We show that the solution of the relaxation actually corresponds to a feasible controller.
\end{enumerate}
\end{remarks}

We would like to emphasize the very important role played by the availability of \emph{full state-feedback}. As a consequence of state-feedback, we have that ${P}_{21}=(sI-A)^{-1}F$. Thus ${P}_{21}$ is  left invertible (though the inverse is improper), and in the (block) incidence algebra. It is this very important feature of ${P}_{21}$ that allows us to use this modified parametrization mentioned \eqref{eq:param} in the preceding paragraph. This parametrization enables us to rewrite the problem in the form \eqref{eq:6a}. This form will turn out to be crucial to our main separability result (Theorem \ref{theorem:2}), which enables us to separate the decentralized problem into a set of decoupled centralized problems.

A main step in our solution strategy will be to reduce the optimal control problem to a set of standard centralized control problems, whose solutions may be obtained by solving standard Riccati equations. The key result about centralized $\mc{H}_{2}$ optimal control is as follows.
\begin{lemma}\label{lemma:centralized}
Consider a system $H$ given by 
$$
H=\ll \ba{cc} H_{11} & H_{12} 
\ea \rr 
= \ll \ba{c|cc} A_{H} & F_{H} & B_{H} \\ \hline
C_H & 0 & D_H 
\ea \rr
$$

along with the following optimal control problem: 
\begin{equation}\label{eq:6}
\ba{rl}
\underset{Q}{\text{minimize}} \text{ \ \ } & \|H_{11}+H_{12}Q \|^{2} \\
\text{subject to \ \ }& Q \text{ stable.}
\ea
\end{equation}
Suppose the pair $(A_H,B_H)$ is stabilizable, $C_H^{T}D_H=0$, and $D_H^{T}D_H\succ0$. 
Then the following Riccati equation has a unique symmetric and positive definite solution $X$:
\begin{equation} \label{eq:12a}
A_H^{T}X+XA_H-XB_H(D_H^{T}D_H)^{-1}B_{H}^{T}X+C_H^{T}C_H=0.
\end{equation}
Let $L$ be obtained from this unique positive definite solution via:
\begin{equation} \label{eq:13a}
L=(D_H^{T}D_H)^{-1}B_H^{T}X.
\end{equation}
Then the optimal solution to \eqref{eq:6} is given by:
\begin{equation}\label{eq:Ricc}
\ba{rl}
Q=\ll \ba{c|c}
A_H-B_HL & F_H \\ \hline -L & 0
\ea \rr.
\ea
\end{equation}
\end{lemma}
(We will often refer to the trio of equations \eqref{eq:12a}, \eqref{eq:13a}, \eqref{eq:Ricc} by $(L,Q)=\text{Ric}(H)$.) 
\begin{proof}
The proof is based on standard techniques and can be argued via a completion-of-squares argument. In particular, it follows from the solution to the standard $\mc{H}_2$ optimal control problem \cite[Theorem 14.7]{ZDG}. Using this theorem, the solution to the $\mc{H}_2$ optimal control problem for the standard problem with the data
$$
G=\ll \ba{cc}
G_{11} & G_{12} \\
I & 0
\ea \rr=
\ll \ba{c|cc}
A_H & F_H & B_H \\ \hline
C_H & 0 & D_H \\
0 & I & 0
\ea \rr
$$
gives the required formula.
\end{proof}

\subsection{Separability of Optimal Control Problem}
We next illustrate the main solution strategy via a simple example.
 Consider the decentralized control problem \eqref{eq:chap4_opt_controll} for the poset in Fig. \ref{fig:2.1}(b). Using the reformulation \eqref{eq:6a} the optimal control problem  \eqref{eq:chap4_opt_controll} may be recast as:
\begin{equation*}
\ba{rl}
\underset{Q}{\text{minimize}} & \left\|P_{11}+P_{12}\ll \ba{ccc} Q_{11} & 0 & 0 \\
Q_{21} & Q_{22} & 0 \\
Q_{31} & 0 & Q_{33}
 \ea \rr 
\right\|^{2} \\
\ea
\end{equation*}
Note that 
$ P_{12}(\down 1)=P_{12}, P_{12}(\down 2)=P_{12}(2)$ (second column of $P_{12}$),  and $P_{12}(\down 3)=P_{12}(3)$.
Similarly $Q^{\down 1}= \ll \ba{ccc} Q_{11}^{T} & Q_{21}^{T} & Q_{31}^{T} \ea \rr^{T}$, $Q^{\down 2}=Q_{22}$, and $Q^{\down 3}=Q_{33}$.
Due to the column-wise separability of the $\mc{H}_{2}$ norm, the problem can be recast as: 

\begin{equation*}
\ba{rl}
\underset{Q}{\text{minimize}} & \left\|P_{11}(1)+P_{12}(\down )\ll \ba{c} Q_{11}  \\
Q_{21}  \\
Q_{31}  \ea \rr  \right\|^{2} +  \left\|P_{11}(2)+P_{12}(\down 2) Q_{22}  \right\|^{2} \\ 
& + \left\|P_{11}(3)+P_{12}( \down 3)  Q_{33}  \right\|^{2} \\
\ea
\end{equation*}\normalsize
Since the sets of variables appearing in each of the three quadratic terms  is different, the problem now may be decoupled into three separate sub-problems, each of which is a standard centralized control problem. For instance, the solution to the second sub-problem can be obtained by noting the realizations of $P_{11}(2)$ and $P_{12}(2)$ and then using \eqref{eq:Ricc}. In this instance, 
$$(G_{22}^{*},Q_{22}^{*})=\text{Ric}\left( \ll \ba{cc} P_{11}(2)& P_{12}(\down 2) \ea \rr \right)=\text{Ric}\left( \ll \ba{c|cc} A_{22} & F_{22} & B_{22} \\ \hline
C(2) & 0 & D(2) 
\ea \rr \right).$$ 
In a similar way, the entire matrix $Q^{*}$ can be obtained, and by design $Q^{*} \in \mc{I}(\mc{P})$ (and is stabilizing). To obtain the optimal $K^{*}$, one can use \eqref{eq:5}. In fact, it is possible to give an explicit state-space formula for $K^{*}$, this is the main content of Theorem~\ref{theorem:3} in the next section.

\section{Main Results}\label{sec:4}
In this section, we present the main results of the paper. The proofs are available in Section \ref{sec:proofs}.
\subsection{Problem Decomposition and Computational Procedure}
\begin{theorem}[Decomposition Theorem]\label{theorem:2}
Let $\mc{P}$ be a poset and $\mc{I}(\mc{P})$ be its incidence algebra. Consider a poset-causal system given by \eqref{eq:1}. The problem  \eqref{eq:6a} is equivalent to the following set of $|P|$ independent decoupled problems:
\begin{equation}\label{eq:reform1}
\begin{aligned}
& \underset{Q^{\down j}}{\text{minimize}} & &   \|P_{11}(j)+P_{12}(\down j)Q^{\down j} \|^{2}  \text{ \ \ \ } \forall j\in P.\\
\end{aligned}
\end{equation}
\end{theorem}
Theorem \ref{theorem:2} is essentially the first step towards a state-space solution. The advantage of this equivalent reformulation of the problem is that we now have $p=|P|$ sub-problems, each over a different set of variables (thus the problem is \emph{decomposed}). Moreover, each sub-problem corresponds to a particular \emph{standard centralized} control problem, and thus the optimal $Q$ in \eqref{eq:chap4_opt_controll} can be computed by simply solving each of these sub-problems.

The subproblems described in \eqref{eq:reform1} have the following interpretation. Once a controller $K$, or equivalently $Q$ is chosen a map $T_{zw}$ from the exogenous inputs $w$ to the outputs $z$ is induced. Let us denote by $T_{zw}(1)$ to be the map from the first input $z_1$ to all the outputs $w$ (this corresponds to the first column of $T_{zw}$). Similarly, the map from $z_i$ to $w$ for $i \in P$ is given by $T_{zw}(i)$. These subproblems correspond to the computation of the \emph{optimal} maps $T_{zw}^{*}(i)
$ for all $i \in P$ from the $i^{th}$ input $z_i$ to the output $w$. The decomposability of the $\mc{H}_2$ norm implies that these maps may be computed separately, and the performance of the overall system is simply the aggregation of these individual maps.


Our next theorem provides an efficient computational technique to obtain the required state-space solution. To obtain the solution, one needs to solve Riccati equations corresponding to the sub-problems we saw in Theorem \ref{theorem:2}. We combine these solutions to form certain simple block matrices, and after simple LFT transformations, one obtains the optimal controller $K^{*}$.

 Before we state the theorem, we introduce some relevant notation. 
\begin{definition}\label{def:Riccati}
We define the operator $\text{Ric}(\down j)$ for $j \in P$ by:
\begin{equation}\label{eq:Ric_short}
\text{Ric}(\down j):=\text{Ric}\left( \ll \ba{c|cc} A(\down j, \down j) & F_{jj} & B(\down j, \down j) \\ \hline 
C(\down j) & 0 & D(\down j)
\ea \rr \right).
\end{equation}
\end{definition}
We define $K(\down j, \down j)$ via $(K(\down j, \down j),Q(j))=\text{Ric}(\down j)$ for $j \in P$.
We introduce two matrices related to the above solution, namely: \begin{align*}
\mathbf{A}&=\text{diag}(A(\down j, \down j)-B(\down j, \down j)K(\down j, \down j)) \\
\mathbf{K}&=\text{diag}(K(\down j, \down j)). 
\end{align*}
We will see later on that $\mathbf{A}$ is the closed-loop state transition matrix under a particular indexing of the states.
We introduce three matrices related to structure of the poset, namely:
\begin{equation} \label{eq:matrix_defs0}
\begin{split}
\Pi_1&=\ll \ba{cccc}  
E_{1} & 0  & \ldots & 0 \\
0 & E_{1} & \ldots & 0 \\
\vdots & \vdots & \ddots & \\
0 & 0 & & E_{1}
\ea \rr, \\
\Pi_2&=\text{diag}\left( \ll \ba{ccc} E_2 & \ldots & E_{|\down j|} \ea \rr \right), \\
R&=\ll \ba{ccc}
E_{\down 1} & \ldots & E_{\down p}
\ea \rr.
\end{split}
\end{equation}
(In $\Pi_1$, the $j^{th}$ diagonal block $E_1$ has $|\down j|$ number of block rows. To be precise, it is a $(\sum_{k \in \down j}n_k)\times n_j$ matrix with the first $n_j \times n_j$ block as the identity and the rest zeroes.) These matrices also have a natural interpretation. In writing the overall states of the closed loop in vector form, we first write the states of subsystem $1$ (i.e. $x_1$) of the plant, then the states of the controller for subsystem $1$ (i.e. $q(1)$), then subsystem $2$ plant states and controller states, and so on. In this indexing, $\Pi_1$ is a \emph{projection operator} that projects onto the coordinates of all the state variables $x_1, \ldots, x_2$.
$\Pi_2$ is simply the matrix that projects onto the orthogonal complement, i.e.  the controller variables $q(1), \ldots, q(p)$. The optimal controller and other related objects can be expressed in terms of the following matrices:
\begin{equation} \label{eq:matrix_defs}
\begin{split}
A_{\Phi}&=\Pi_2^{T}\mathbf{A}\Pi_2, \\
B_\Phi&=\Pi_2^{T}\mathbf{A}\Pi_1, \\
C_\Phi&=R\Pi_2, \\
C_Q&=-R\mathbf{K}.
\end{split}
\end{equation}
We illustrate this notation further by means of a numerical example in Appendix \ref{sec:num}.


\begin{theorem}[Computation of Optimal Controller] \label{theorem:3}
Consider the poset-causal system of the form \eqref{eq:1}, with $(A_{ii}, B_{ii})$ stabilizable for all $i \in P$. Consider the following Riccati equations: 
$$
\left(K(\down j, \down j),Q(j) \right)=\text{Ric}\left(\down j\right) \; \; \forall j \in P.
$$
Then the optimal solution to the problem \eqref{eq:chap4_opt_controll} is given by the controller:
\begin{equation}\label{eq:thm2}
K^{*}= \ll \ba{c|c}
A_{\Phi}-B_{\Phi}C_{\Phi} & B_{\Phi} \\ \hline
C_Q(\Pi_2-\Pi_1C_{\Phi}) &C_Q\Pi_1
\ea \rr.
\end{equation}
Moreover, the controller $K^{*} \in \mc{I}(\mc{P})$ and is internally stabilizing.
\end{theorem}

Recall that $n_i$ denotes the degree of the $i^{th}$ sub-system in \eqref{eq:2}. Let $n_{\max}=\max_i n_i$ be the largest degree of the sub-systems. Let $n(\down \down i)=\sum_{j \in \down \down i} n_j$. Let $\sigma_{\mc{P}}=\sum_{j \in P} |\down \down j|$ (note that this is a purely combinatorial quantity, dependent only on the poset).
As we mentioned in the introduction, one of the advantages of state-space techniques is that they provide graceful degree bounds for the optimal controller. As a consequence of Theorem~\ref{theorem:3} we have the following:
\begin{cor}[Degree Bounds]\label{cor:degbounds}
The degree $d_{K^{*}}$ of the overall optimal controller is bounded above by 
$$
d_{K^{*}} \leq \sum_{j \in P} n(\down \down j).
$$
In particular, 
$
d_{K^{*}}  \leq \sigma_{\mc{P}}n_{\max}.
$ 
Moreover, the degree of the controller implemented by subsystem $j$ is bounded above by $n(\down \down j)$.
\end{cor}

\subsection{Structure of the Optimal Controller} \label{sec:structure_optimal_controller}
Having established the computational aspects, we now turn to some structural aspects of the optimal controller. We first introduce a pair of very important objects $(\Phi, \Gamma)$, called the \emph{propagation filter} and the \emph{differential filter}, respectively.  Define the block $p \times p$ transfer function matrices $(\Phi, \Gamma)$ via:
\begin{align} \label{eq:phi_gamma}
\Phi=\ll  \ba{c|c} 
\text A_{\Phi} & B_{\Phi} \\ \hline 
C_{\Phi} & I
\ea \rr & \text{ \ \ \ \ \ \ \ \ \ }
\Gamma=\ll  \ba{c|c} 
\text A_{\Phi} - B_{\Phi}C_{\Phi} & B_{\Phi} \\ \hline 
-C_{\Phi} & I
\ea \rr. 
\end{align}
Note that both $\Phi$ and $\Gamma$ are invertible (since their ``$D$'' matrices are
equal to $I$), and in fact, they are inverses of each other, i.e., $\Gamma \Phi =
\Phi \Gamma = I$.   We sometimes denote the entries $\Phi_{ij}=\Phi_{i \leftarrow j}$ and similarly  $\Gamma_{ij}=\Gamma_{i \leftarrow j}$ to emphasize a certain interpretation of these quantities. We note that $\Phi_{i \leftarrow i}= \Gamma_{i \leftarrow i}=I$ (this can be seen from the fact that the corresponding entries in the ``$C$'' matrices of the transfer functions is zero). We show (Lemma \ref{lemma:inc_alg} in Section \ref{sec:proofs}), that $\Phi, \Gamma \in \mc{I}(\mc{P})$. Moreover, the fact that $\Phi^{-1}=\Gamma$ in conjunction with Corollary \ref{cor:inc} gives the following expression:
\begin{equation} \label{eq:19}
\Gamma_{i \leftarrow j}=\sum_{p_{ij} \in [j \rightarrow i]} \prod_{\{l,k \} \in p_{ij}}(-\Phi_{l \leftarrow k}).
\end{equation}

We will show that $\Phi_{l \leftarrow k}$, in fact, corresponds to a \emph{specific filter} that propagates local signals upstream. For example, in Fig. \ref{fig:2.1}(a), if $x_1$ is the state at subsystem $1$, $\Phi_{21}x_1$ is the prediction of state $x_2$ at subsystem $1$. On the other hand, $\Gamma$  has an interesting dual interpretation. As one proceeds ``upstream'' through the poset, more information is available, and consequently the prediction of the global state becomes more accurate. The transfer function $\Gamma$ plays the role of computing the differential improvement in the prediction of the global state. For this reason, we call it the \emph{differential filter}. Interestingly, it is intimately related to the notion of \emph{M\"{obius} inversion} on a poset, a generalization of differentiation to posets.
%
We briefly discuss these ideas in the ensuing discussion. Before stating the next theorem, we introduce the transfer function matrix $K_{\Phi}$, which is defined column-wise via:
\begin{equation} \label{eq:Kphi}
K_{\Phi}(j)=\hat{K}(\down j, \down j)\Phi(j).
\end{equation}

\begin{theorem}[Structure of Optimal Controller]\label{theorem:4}
The optimal controller \eqref{eq:thm2} is of the form:
\begin{align*}
u(t)&=-K_{\Phi}\Gamma x(t) \\
&=-\sum_{j \in P} \hat{K}(\down j, \down j) \Phi(j) (\Gamma x)_j(t).
\end{align*}
\end{theorem}
\begin{remark} 

Let us denote the vector $e(j)=\Phi(j)(\Gamma x)_j$. We will interpret $e(j)$ as the differential improvement in the prediction of the global state $x$ at subsystem $j$. Denoting $\hat{K}(\down j, \down j)$ by $K_j$, note that the control law takes the form $u(t)=\sum_{j \in P} K_j e(j)$. This structural form suggests that the controller uses the differential improvement of the global state at the different subsystems as the atoms of local control laws, and that the overall control law is an aggregation of these local control laws.
\end{remark}
\subsection{Interpretation of $\Phi$ and $\Gamma$}\label{sec:phi_gamma_interpretation}
In this section we explain the role of the propagation filter $\Phi$ and the differential filter $\Gamma$.
Due to the information constraints in the problem, at subsystem $j$ only states in $\up j$ are available, states of other subsystems are unavailable. A reasonable architecture for the controller would involve \emph{predicting} the unknown states at subsystem $j$ from the available information. This is illustrated by the following example.

\begin{example}
Consider the system shown in Fig. \ref{fig:chap4_partial_pred} with dynamics
$$
\left[ \begin{array}{c} \dot{x}_1  \\ \dot{x}_2\\ \dot{x}_3 \end{array} \right]= \left[ \begin{array}{ccc} A_{11} & 0 & 0 \\ 0 & A_{22} & 0 \\ A_{31} & A_{32} & A_{33} \end{array}\right]
\left[ \begin{array}{c} x_1  \\ x_2\\ x_3 \end{array} \right] + \left[ \begin{array}{ccc} B_{11} & 0 & 0 \\ 0 & B_{22} & 0 \\ B_{31} & B_{32} & B_{33} \end{array}\right]\left[ \begin{array}{c} u_1  \\ u_2\\ u_3 \end{array} \right].
$$

\begin{figure}[htbp]
  \begin{center}
  \includegraphics[scale=0.5]{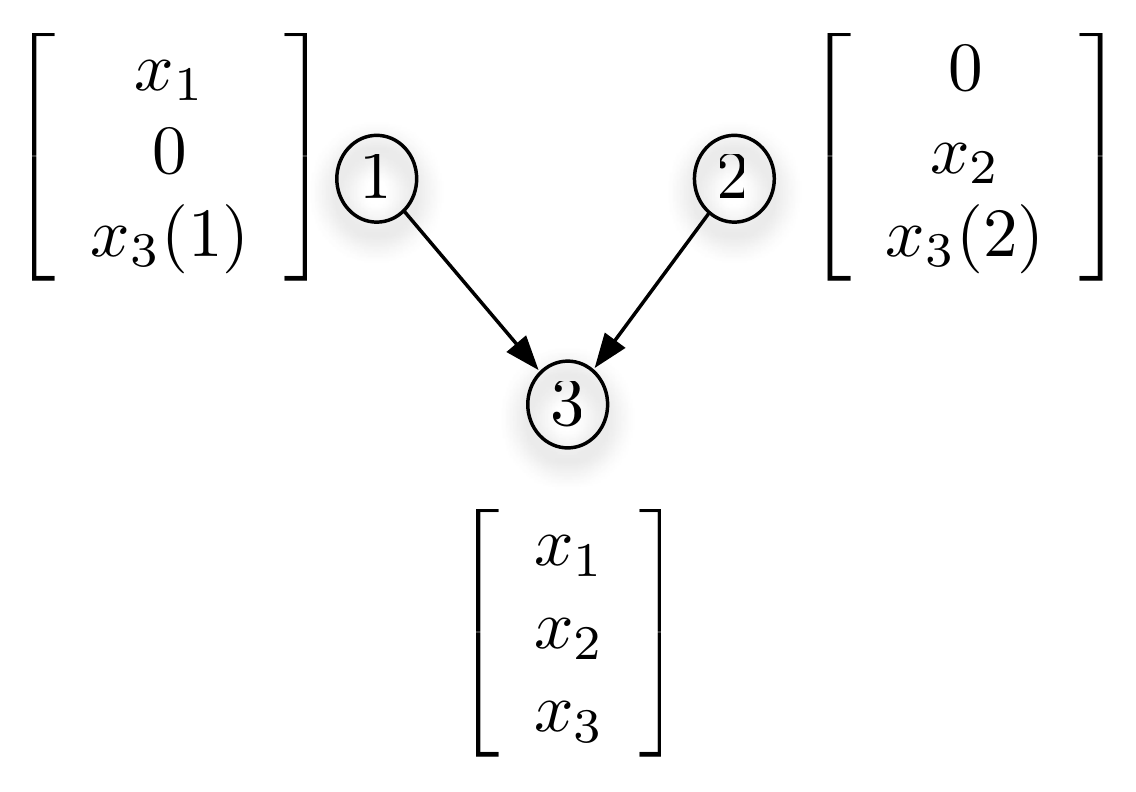}
  \end{center}
  \caption{Local state information at the different subsystems. The quantities $x_3(1)$ and $x_3(2)$ are partial state predictions of $x_3$.}
  \label{fig:chap4_partial_pred}
  \end{figure}
  
Note that subsystem $1$ has no information about the state of subsystem $2$. Moreover, the state $x_1$ or input $u_1$ do not affect the dynamics of $2$ (their respective dynamics are uncoupled). Hence the only sensible prediction of $x_2$ at subsystem $1$ (which we denote by $x_2(1)$) is $x_2(1)=0$. Subsystem $1$ also does not have access to the state $x_3$. However, it can predict $x_3$ based on the influence that the state $x_1$ has on $x_3$. (Note that both the states $x_1$, $x_2$ and inputs $u_1$, $u_2$ affect $x_3$ and $u_3$.) Let us denote $x_3(1)$ to be the prediction of state $x_3$ at subsystem $1$. Since $x_2$ and $u_2$ are unknown, the state $x_3(1)$ is a \emph{partial} prediction of $x_3$ (i.e. $x_3(1)$ is the prediction of the \emph{component of} $x_3$ that is affected by subsystem $1$). Similarly, subsystem $2$ maintains a prediction of $x_3$ denoted by $x_3(2)$, which is also a partial prediction of $x_3$. Each subsystem thus maintains (possibly partial) predictions of unknown downstream states, as shown in Fig. \ref{fig:chap4_partial_pred}.
\end{example}

In this paper we will not discuss how the state predictions are computed, a detailed discussion of the same is available in \cite[Chapter 5]{pari_thesis}, \cite{poset_struct}. However, we mention that $\Gamma$ has an interesting related role. At subsystem $j$ the true state $x_j$ becomes available for the first time (with respect to the subposet $\up j$). The quantity $q_j(j):=(\Gamma x)_j$ measures the \emph{differential improvement} in the knowledge of state $x_j$, i.e. the difference between the true state $x_j$ and its best prediction from upstream information. Similarly, it is possible do define $q_j(i)$, the differential improvement in the prediction of state $j$ at subsystem $i$. We let
$$
q(i)=[q_j(i)]_ {j \in \down \down i},
$$
so that $q(i)$ corresponds to the differential improvement in state predictions at the $i^{th}$ subsystem. This $q(i)$ is a vector of length $|\down \down i|$ and its components $q_j(i)$ are differential improvements of $x_j$ for $j \in \down \down i$ at subsystem $i$.

We next examine the role of $\Phi$. Consider a system of the form:
$$
\dot{r}(t)=Hr(t),
$$
where $r(t)\in \mathbb{R}^{n}$. Given $r_1(t)$ it is possible to compute $r_2(t), \ldots, r_n(t)$ by \emph{propagation} by noting that $(sI-H)r(z)=0$. Rewriting these equations, we obtain that
$$
\ll \ba{cc} sI-H_{11} & -H_{12} \\ -H_{21} & sI-H_{22} \ea \rr \ll \ba{c} r_1 \\ \vdots \\r_n \ea \rr=0
$$
to obtain
\begin{equation}\label{eq:phi1}
\ll \ba{c}  r_2 \\ \vdots \\ r_n \ea \rr(s)=(sI-H_{22})^{-1}H_{21}r_1(s),
\end{equation}
where $H_{22}=E_{\{2, \ldots, n\} }^{T}HE_{\{2, \ldots, n\} }$ and $H_{21}=E_{\{2, \ldots, n\} }^{T}HE_{1}$.
The map $\Phi_{H}=(sI-H_{22})^{-1}H_{21}$ from $r_1$ to $r_2, \ldots, r_n$ is simply a propagation of the ``upstream'' states based on $r_1$.

It is possible to show \cite[Chapter 5]{pari_thesis}, \cite{poset_struct} that the \emph{differential improvements} in the local state predictions $q(i)$ obey a decoupled relationship as a consequence of a separation principle. (Thus $r$ in \eqref{eq:phi1} corresponds to the differential improvement $q(i)$ for some subsystem $i$).
As we already mentioned $(\Gamma x)_i$ is the differential improvement in $x_i$ at subsystem $i$. Since $\Phi$ plays the role of propagating decoupled local signals, it follows that $\Phi_{ji}(\Gamma x)_i$ is the differential improvement in the prediction of the state $x_j$ for $j \in \down \down i$ at subsystem $i$. More precisely,
\begin{equation} \label{eq:q_phi}
q_j(i)=\Phi_{ji}(\Gamma x)_i, 
\end{equation}
and this can be written compactly as $q(i)=\Phi(i) (\Gamma x)_i$, where the column vector $q(i)$ is the differential improvement in the global state at subsystem $i$. As an aside, we mention that the states of the optimal controller correspond precisely to these differential improvements $q_j(i)$.
\subsection{Structure of the Optimal Controller}
Using Theorem \ref{theorem:4}, the optimal control law can be expressed as:
\begin{equation} \label{eq:interpret5}
u=\sum_{i \in P} \hat{K}(\down i, \down i)\Phi(i) (\Gamma x)_i. 
\end{equation}
As explained above, $\Phi(i) (\Gamma x)_i$ is a vector containing the differential improvement in the prediction of the global state at subsystem $i$. Each term $\hat{K}(\down i, \down i)\Phi(i) (\Gamma x)_i$ may be viewed as a local control law acting on the local differential improvement in the predicted state. The overall control law has the elegant interpretation of being an aggregation of these local control laws.


\begin{example}
Let us consider the poset from Fig. \ref{fig:2.1}(d), and examine the structure of the controller. (For simplicity, we let $K_j=\hat{K}(\down j, \down j)$, the gains obtained by solving the Riccati equations). The control law may be decomposed into local controllers as:
\begin{align*}
u&=K_1\ll \ba{c} I  \\ \Phi_{21} \\ \Phi_{31} \\ \Phi_{41} \ea \rr x_1+
K_2\ll \ba{c} 0  \\ I \\ 0 \\ \Phi_{42} \ea \rr (\Gamma x)_2+
K_3\ll \ba{c} 0  \\ 0 \\ I \\ \Phi_{43} \ea \rr (\Gamma x)_3+
K_4 \ll \ba{c} 0  \\ 0 \\ 0 \\ I \ea \rr (\Gamma x)_4 \\
&=K_1\ll \ba{c} x_1 \\q_2(1) \\ q_3(1)\\ q_4(1) \ea \rr +
K_2\ll \ba{c} 0  \\ x_2-q_2(1) \\ 0 \\ q_4(2) \ea \rr +
K_3\ll \ba{c} 0  \\ 0 \\ x_3-q_3(1) \\ q_4(3) \ea \rr +
K_4 \ll \ba{c} 0  \\ 0 \\ 0 \\ (x_4-q_4(1))-q_4(2)-q_4(3) \ea \rr. \\
\end{align*}
Each term in the above expression has the natural interpretation of being a local control signal corresponding to differential improvement in predicted states, and the final controller can be viewed as an aggregation of these. 

Note that zeros in the above expression imply no improvement on the local state. For example, at subsystem $2$ there is no improvement in the predicted value of $x_3$ because the state $x_2$ does not affect subsystem $3$ due to the poset-causal structure. There is no improvement in the predicted value of state $x_3$ at subsystem $4$ either,  because the best available prediction of $x_3$ from downstream information $\up \up 4$ is $x_3$ itself. While this interpretation has been stated informally here, it has been made precise in \cite[Chapter 5]{pari_thesis}, \cite{poset_struct}.
\end{example}

\section{Discussion and Examples} \label{sec:5}


\subsection{The Nested Case}
Consider the poset on two elements $\mc{P}=\left( \l 1,2  \r , \preceq \right)$ with the only order relation being $1\preceq 2$ (Fig. \ref{fig:2.1}(a)). This is the poset corresponding to the communication structure in the ``Two-Player Problem'' considered in \cite{Swigart}. We show that their results are a specialization of our general results in Section \ref{sec:4} restricted to this particular poset.

We begin by noting that from the problem of designing a nested controller (again we assume $F=I$ for simplicity) can be recast as:
\begin{equation*}
\ba{rl}
\underset{Q}{\text{minimize}} & \left\|P_{11}+P_{12}\ll \ba{cc} Q_{11} & 0 \\
Q_{21} & Q_{22} \\
 \ea \rr \right\|^{2} \\
\ea
\end{equation*}

By Theorem \ref{theorem:2} this problem can be recast as:
\begin{equation*}
\ba{rl}
\underset{Q}{\text{minimize}} & \left\|P_{11}^{1}+P_{12}\ll \ba{c} Q_{11}  \\
Q_{21}  \ea \rr \right\|^{2} +  \left\|P_{11}^{2}+P_{12}^{2} Q_{22}  \right\|^{2} \\ 
\ea
\end{equation*}\normalsize

We wish to compare this to the results obtained in \cite{Swigart}. It is possible to obtain precisely this same decomposition in the finite time horizon where the $\mc{H}_2$ norm can be replaced by the Frobenius norm and separability can be used to decompose the problem. For each of the sub-problems, the corresponding optimality conditions may be written (since they correspond to simple constrained-least squares problems). These optimality conditions correspond exactly to the decomposition of optimality conditions they obtain (the crucial Lemma 3 in their paper). We point out that the decomposition is a simple consequence of the separability of the Frobenius norm. 

Let us now examine the structure of the optimal controller via Theorem \ref{theorem:4}. Note that $\down 1= \l 1,2 \r$ and $\down 2 =\l 2 \r$. Based on Theorem \ref{theorem:3}, we are required to solve $(K,Q(1))=\text{Ric}(\down 1)$, and $(J,Q(2))=\text{Ric}(\down 2)$.
Noting that in this example $\Gamma_{2 \leftarrow 1}=-\Phi_{2 \leftarrow 1}$, a straightforward application of Theorem \ref{theorem:4} yields the following:
\begin{align*}
u_{1}(t)&=-(K_{11}+K_{12}\Phi_{2 \leftarrow 1})x_{1}(t) \\
u_{2}(t)&=-(K_{21}+K_{22}\Phi_{2 \leftarrow 1})x_{1}(t)-J(x_{2}(t)-\Phi_{21}x_1(t)), \\
\end{align*}
which is precisely the structure of the optimal controller given in \cite{Swigart,Swigart2}. It is possible to show (as Swigart et. al indeed do in \cite{Swigart}) that $\Phi_{2 \leftarrow 1}$ is an predictor of $x_{2}$ based on $x_{1}$. Thus the controller for $u_1$ predicts the state of $x_2$ from $x_1$, uses the estimate as a surrogate for the actual state, and uses the gain $K_{21}$ in the feedback loop. The controller for $u_2$ (perhaps somewhat surprisingly) also estimates the state $x_2$ based on $x_1$ using $\hat{x}_{2}=\Phi_{21}x_1$ (this can be viewed as a ``simulation'' of the controller for $u_1$). The prediction error for state $2$ is then given by $e_2:=x_2-\hat{x}_2=x_2-\Phi_{21}x_1$. The control law for $u_2$ may be rewritten as 
$$
u_2=-(K_{21}x_1+K_{22}\hat{x}_2+Je_2).
$$
Thus this controller uses predictions of $x_2$ based on $x_1$ along with prediction errors in the feedback loop. We will see in a later example, that this prediction of states higher up in the poset is prevalent in such poset-causal systems, which results in somewhat larger order controllers. 

Analogous to the results in \cite{Swigart}, it is possible to derive the results in this paper for the finite time horizon case (this is a special case corresponding to FIR plants in our setup). We do not devote attention to the finite time horizon case in this paper, but just mention that similar results follow in a straightforward manner. 


\subsection{Discussion Regarding Computational Complexity}
Note that the main computational step in the procedure presented in Theorem \ref{theorem:3} is the solution of the $p$ sub-problems. The $j^{th}$ sub-problem requires the solution of a Riccati equation of size at most $|\down j| n_{\max}=O(p)$ (when the degree $n_{\max}$ is fixed). Assuming the complexity of solving a Riccati equation using linear algebraic techniques is $O(p^4)$ \cite{Carroll} the complexity of solving $p$ of them is at most $O(p^{5})$. We wish to compare this with the only other known state-space technique that works on all poset-causal systems, namely the results of Rotkowitz and Lall \cite{quadinv}. In this paper, they transform the problem to a standard centralized problem using Kronecker products. In the final computational step, one would be required to solve a single large Riccati equation of size $O(p^2)$, resulting in a computational complexity of $O(p^{8})$. 
\subsection{Discussion Regarding Degree Bounds}
It is insightful to study the asymptotics of the degree bounds in the setting where the sub-systems have fixed degree and the number of sub-systems $p$ grows. As an immediate consequence of the corollary, the degree of the optimal controller (assuming that the degree of the sub-systems $n_{\max}$ is fixed) is at most $O(p^{2})$ (since $n(\down j) \leq p$). In fact, the asymptotic behaviour of the degree can be sub-quadratic. Consider a poset $(\l1, \ldots, p \r, \preceq)$ with the only order relations being $1 \preceq i$ for all $i$. Here $|\down 1| =p$, and $|\down i| = 1$ for all $i \neq 1$. Hence, $\sum_{j}|\down j| - p \leq p$, and thus $d^{*} \leq sn_{\max}$. In this sense, \emph{the degree of the optimal controller is governed by the poset parameter $\sigma_{\mc{P}}$}.

\section{Proofs of the Main Results} \label{sec:proofs}
\begin{proof}[Proof of Theorem \ref{theorem:stab}]
Note that one direction is trivial. Indeed if the $(A_{ii}, B_{ii})$ are stabilizable, one can pick a diagonal controller with diagonal elements $K_{ii}$ such that $A_{ii}+ B_{ii}K_{ii}$ is stable for all  $i \in P$. This constitutes a stabilizing controller.

For the other direction let 
$$
K=\ll \ba{c|c}
A_K & B_K \\ \hline
C_K & D_K
\ea \rr
$$
be a poset-causal controller for the system. We will first show that without loss of generality, we can assume that $A_K, B_K, C_K, D_K$ are block lower triangular (so that $K$ has a realization where all matrices are block lower triangular).

First, note that since $K \in \mc{I}(\mc{P})$, $D_K \in \mc{I}(\mc{P})$. Recall, that we assumed throughout that the indices of the matrices in the incidence algebra are labeled so that they are consistent with a linear extension of the poset, so that $D_K$ is lower triangular. Note that the controller $K$ is a block $p \times p$ transfer function matrix which has a realization of the form:

$$
K=\ll \ba{c|ccc}
A_K & B_K(1) & 
\ldots & B_K(p) \\ \hline
C_K(1) & D_K(1,1) & \ldots & D_K(1,p) \\
\vdots & \vdots & \ddots & \\
C_K(p) & D_K(p,1) & \ldots & D_K(p,p)
\ea \rr
$$
Since the controller $K \in \mc{I}(\mc{P})$, we have that $K_{pj}=0$ for all $j \neq p$ (recall that $p$ is the cardinality of the poset). This vector of transfer functions (given by the last column of $K$ with the $(p,p)$ entry deleted) is given by the realization:
$$
\bar{K}_p:=\ll \ba{c} C_K(1) \\ \vdots \\ C_{K}(p-1) \ea \rr (sI-A_K)^{-1}B_{K}(p) +\ll \ba{c} D_K(1,p) \\ \vdots \\ D_{K}(p-1,p) \ea \rr=0.
$$

Since this transfer function is zero, in addition to $D_{K}(j,p)=0$ for all $j=1, \ldots, p-1$, it must also be the case that the controllable subspace of $(A_K, B_K(p))$ is contained within the unobservable subspace of $\left( \left[ \ba{ccc} C_K(1)^{T} & \ldots & C_K(p-1)^{T} \ea \right]^{T}, A_K \right)$. 
By the Kalman decomposition theorem \cite[pp. 247]{desoer}, there is a realization of this system of the form:
$$
\bar{K}_p=\ll \ba{c|c} \bar{A} & \bar{B} \\ \hline \bar{C} & \bar{D} \ea \rr,
$$
where $(\bar{A},\bar{B},\bar{C},\bar{D})$ are of the form:
\begin{equation} \begin{split}\label{eq:block_structure}
\bar{A}&=\ll \ba{ccc} A_{11} & 0 & 0 \\ 
A_{21}& A_{22} & 0 \\
A_{21} & A_{32} & A_{33} \ea \rr \\
\bar{B}&=\ll \ba{c}0 \\ 0 \\ B_3  \ea \rr \\
\bar{C}&=\ll \ba{ccc} C_1 & 0 & 0 \ea \rr \\
\bar{D}&=0.
\end{split}
\end{equation}
 As an aside, we remind the reader that this decomposition has a natural interpretation. For example, the subsystem 
$$\left( \ll \ba{cc} A_{11} & 0 \\ A_{21} & A_{22} \ea \rr, \ll \ba{c}0 \\ 0 \ea \rr,\ll \ba{cc} C_1 & 0 \ea \rr \right)
$$
corresponds to the observable subspace, where the system is uncontrollable, etc. 
(The usual Kalman decomposition as stated in standard control texts is a block $4 \times 4$ decomposition of the state-transition matrix. Here we have a smaller block $3 \times 3$ decomposition because of the collapse of the subspace where the system is required to be both controllable and observable).

Thus this decomposition allows us to infer the specific block structure \eqref{eq:block_structure} on the matrices $(\bar{A},\bar{B},\bar{C},\bar{D})$.
As a result of this block structure, there is a realization of the overall controller  $(A_K, B_K, C_K, D_K)$, where all the matrices have the block structure
$$
\ll \ba{cccc}
M_{1,1} & \ldots & M_{1,p-1} & 0 \\
\vdots & \ddots & &\vdots \\
M_{p-1,1} & \ldots & M_{p-1,p-1} & 0 \\
M_{p,1} & \ldots & M_{p,p-1} & M_{p,p} 
\ea \rr.
$$
One can now repeat this argument for the upper $(p-1)\times (p-1)$ sub-matrix of $K$.
By repeating this argument for first $p-1, p-2, \ldots, 1$ we obtain a realization of $K$ where all four matrices are block lower triangular.

Note that given the controller $K$ (henceforth assumed to have a lower triangular realization), the closed loop matrix $A_{cl}$ is given by
$$
A_{cl}=\ll \ba{cc}
A+BD_K& BC_K \\
B_K & A_K
\ea \rr.
$$
By assumption the (open loop) system is poset-causal, hence $A$ and $B$ are block lower triangular. As a result, each of the blocks $A+BD_K$, $BC_K$, $B_K$, $A_K$ are block lower triangular. A straightforward permutation of the rows and columns enables us to put $A_{cl}$ into block lower triangular form where the diagonal blocks of the matrix are given by
\begin{equation} \label{eq:local_closed_loop}
\ll \ba{cc}
A_{jj}+B_{jj}D_{K_{jj}} & B_{jj}C_{K_{jj}} \\
B_{K_{jj}} & A_{K_{jj}}
\ea \rr. 
\end{equation}
Note that the eigenvalues of this lower triangular matrix (and thus of $A_{cl}$, since permutations of rows and columns are spectrum-preserving) are given by the eigenvalues of the diagonal blocks. The matrix $A_{cl}$ is stable if and only if all its eigenvalues are within the unit disk in the complex plane, i.e. the above blocks are stable for each $j\in P$. 
Note that \eqref{eq:local_closed_loop} is obtained as the closed-loop matrix precisely by the interconnection of $$\ll \ba{c|c} A_{jj}& B_{jj} \\ \hline I & 0 \ea \rr$$ with the controller $$\ll \ba{c|c} A_{K_{jj}}& B_{K_{jj}} \\ \hline C_{K_{jj}} & D_{K_{jj}} \ea \rr.$$ 
Hence, \eqref{eq:local_closed_loop} (and thus the overall closed loop) is stable if and only if $(A_{jj}, B_{jj})$ are stabilizable for all $j \in P$, and $(A_{K_{jj}}, B_{K_{jj}}, C_{K_{jj}}, D_{K_{jj}})$ are chosen to stabilize the pair.
\end{proof}

\begin{proof}[Proof of Theorem 2]
If $G=\ll G_1, \ldots, G_k \rr$ is a matrix with $G_i$ as its columns, then 
\begin{equation*} 
\|G \|_F^2=\sum_{i=1}^{k}\|G_i \|_F^2,
\end{equation*}
where $\| \cdot \|_F$ denotes the Frobenius norm.
This separability property of the Frobenius norm immediately implies the following separability property for the $\mc{H}_2$ norm:
If $H=\ll H_1, \ldots H_k \rr$ is a transfer function matrix with $H_i$ as its columns, then 
\begin{align*} 
\|H \|^2&=\int_{j\mathbb{R}} \|H(\omega)\|_{F}^{2}d\omega 
=\sum_{i=1}^{k}\int_{j\mathbb{R}} \|H_i(\omega)\|_{F}^{2}d\omega
=\sum_{i=1}^{k}\|H_i \|^2,
\end{align*}
(In the above $j\mathbb{R}$ denotes the imaginary axis in the complex plane).
The separability property of the $\mc{H}_2$ norm can be used to simplify \eqref{eq:6}. 
Recall that $P_{11}(j),Q(j)$ denote the $j^{\text{th}}$ columns of $P_{11}$ and $Q$ respectively. Using 
the separability we can rewrite \eqref{eq:6} as
\begin{equation}\label{eq:7b}
\ba{rl}
\underset{Q}{\text{minimize}} & \sum_{j \in P} \|P_{11}(j)+P_{12}Q(j) \|^{2} \\
\text{subject to \ \ } & Q(j) \in \mc{I}(\mc{P})^{j}\\
\ea
\end{equation}
The formulation in \eqref{eq:7b} can be further simplified by noting that for $Q^j \in \mc{I}(\mc{P})^{j}$, 
\begin{equation}\label{eq:8}
P_{12}Q(j)=P_{12}({\down j})Q^{\down j}.
\end{equation}

The advantage of the representation \eqref{eq:8} is that, in the right hand side the variable $Q^{\down j}$ is \emph{unconstrained}. Using this we may reformulate \eqref{eq:7b} as:
\begin{equation}\label{eq:9}
\ba{rl}
\underset{Q^{\down j}}{\text{minimize}} & \sum_{j \in P} \|P_{11}(j)+P_{12}(\down j)Q^{\down j} \|^{2} \\
\ea
\end{equation}
Since the variables in the $Q^{ \down j}$ are distinct for different $j$, this problem can be separated into $p$ sub-problems as follows: 
\begin{equation}\label{eq:10}
\ba{rl}
\underset{Q^{\down j}}{\text{minimize}} & \|P_{11}(j)+P_{12}(\down j)Q^{\down j}\|^{2} \\
&\text{ for all } j\in P.
\ea
\end{equation}
\end{proof}

Note that each sub-problem is a \emph{standard} $\mc{H}_{2}$ optimal centralized control problem, and can be solved using canonical procedures. Once the optimal $Q$ is obtained by solving these sub-problems, the optimal controller may be synthesized using \eqref{eq:5}. The following lemma describes the solutions to the individual sub-problems \eqref{eq:reform1} in Theorem \ref{theorem:2}.
\begin{lemma} \label{lemma:4} Let
$(A,B,C,D)$ be as given in \eqref{eq:2} with $A, B$ in the block incidence algebra $\mc{I}(\mc{P})$. Let
\begin{equation} \label{eq:15.1} (K(\down j, \down j),Q(j))=\text{Ric}(\down j).\end{equation}
 Then the optimal solution of each sub-problem \eqref{eq:reform1} is given by:
 \begin{equation} \label{eq:15}
  (Q^{\down j})^{*}=
 \ll \ba{c|c}
 A(\down j, \down j)-B(\down j, \down j)K(\down j, \down j) & E_1F_{jj} \\ \hline
 -K(\down j, \down j) & 0
 \ea \rr.
 \end{equation}
 \end{lemma} 
 
 (We remind the reader that in the above $E_1$ is the block $| \down j | \times 1 $ matrix which picks out the first column corresponding 
 of the block $|\down j| \times |\down j |$ matrix before it.) \\
\begin{proof}
The proof follows directly from Lemma \ref{lemma:centralized} by choosing 
$$
H=\left[ \ba{cc} P_{11}(j) & P_{12}(\down j)\ea\right]=\ll \ba{c|cc} A(\down j, \down j) & E_1 F_{jj} & B(\down j, \down j) \\ \hline
C(\down j) & 0 & D(\down j) 
\ea \rr.
$$
\end{proof}

\begin{lemma}\label{lemma:A4}
The optimal solution to \eqref{eq:6a} is given by
\begin{equation}\label{eq:A2}
Q^{*}=\left[ \ba{c|c} \mathbf{A} & \Pi_1F\\ \hline
C_Q & 0 \ea \right].
\end{equation}
\end{lemma}
\begin{proof}
We note that Lemma \ref{lemma:4} gives an expression for the individual columns of $Q^{*}$. Using Lemma \ref{lemma:4} and the LFT formula for column concatenation:
$$
\ll \ba{cc} G_1& G_2 \ea \rr = \ll  \ba{cc|cc}
A_1 & 0 & B_1 & 0 \\
0 &A_2 & 0 & B_2 \\ \hline
C_1 & C_2 & D_1 & D_2 
\ea
\rr, 
$$
we obtain the required expression. 
\end{proof}

\begin{lemma}\label{lemma:inc_alg}
The transfer function matrices $\Phi$, $\Gamma$ and $K_{\Phi}$ as given in \eqref{eq:phi_gamma}, \eqref{eq:Kphi}, are in the incidence algebra $\mc{I}(\mc{P})$.
\end{lemma}
\begin{proof}
Let us define block $p \times 1$ transfer functions as follows:
\begin{equation}\label{eq:lem_inc}
\begin{split}
\Phi(j)&= \ll \ba{c|c}
A_{\Phi}(j) & B_{\Phi}(j) \\ \hline
E_{\down \down j} & I
\ea \rr \\
K_{\Phi}(j)&= \ll \ba{c|c}
A_{\Phi}(j) & B_{\Phi}(j) \\ \hline
-\hat{K}(\down j, \down j)E_{\down \down j} & -\hat{K}(\down j, \down j)E_j
\ea \rr 
\end{split}
\end{equation}
Note that $K_{\Phi}(j)=\hat{K}(\down j, \down j) \Phi(j)$. Also, note that if $i$ is such that $j \npreceq i$ then the $i^{th}$ entry of $\Phi(j)$ is zero since the corresponding row of $E_{\down \down j}$ is zero. By similar reasoning, $K_{\Phi}(j)$ also has this property. Thus, when we construct the matrices 
\begin{align*}
 \Phi &= \ll \ba{ccc} \Phi(1) & \ldots & \Phi(p) \ea \rr \\  
 K_{\Phi}&=\ll \ba{ccc} K_\Phi(1) & \ldots & K_\Phi(p) \ea \rr 
 \end{align*}
by column concatenation, we see that both $\Phi \in \mc{I}(\mc{P}) $ and $K_{\Phi} \in \mc{I}_{K}(\mc{P})$. Since $\Gamma=\Phi^{-1}$, we have $\Gamma \in \mc{I}(\mc{P})$.
\end{proof}

%
%

\begin{lemma}\label{lemma:diagstab}
The matrix $\mathbf{A}$ is stable.
\end{lemma}
\begin{proof}
Recall that $\mathbf{A}=\text{diag}(A(\down j, \down j)-B(\down j, \down j)K(\down j, \down j))$.
Since $A(\down j, \down j)$ and $B(\down j, \down j)$ are lower triangular with $A_{kk}$, $B_{kk}$, $k\in \down j$ along the diagonals respectively, we see that the pair $(A(\down j, \down j), B(\down j, \down j))$ is stabilizable by Assumption 1 (simply picking a diagonal $K$ which stabilizes the diagonal terms would suffice to stabilize $(A(\down j, \down j), B(\down j, \down j))$). Hence, there exists a stabilizing solution to $\text{Ric}(\down j)$ and the corresponding controller $K(\down j, \down j)$ is stabilizing. Thus $A(\down j, \down j)-B(\down j, \down j)K(\down j, \down j))$ is stable, and thus so is $\mathbf{A}$.
\end{proof}

Given transfer functions $M$ and $K$, their feedback interconnection is usually described through a linear fractional transformation of the form $
f(M,K)=M_{11}+M_{12}K(I-M_{22}K)^{-1}M_{21}$. State space formulae for this interconnection are standard  \cite[pp. 179]{Zhou} and will be useful for evaluating several quantities in what follows.

\begin{lemma}\label{lemma:LFT}
Given transfer function matrices $M$ and $K$ with realizations 
\begin{align*}
M=\ll \ba{c|cc}
A & B_1 & B_2 \\ \hline
C_1 & D_{11} & D_{12} \\
C_2 & D_{21} & 0
\ea \rr \text{, \ } &
K=\ll \ba{c|c}
A_K & B_K \\ \hline
C_K & D_K
\ea \rr,
\end{align*}
the Linear Fractional Transformation (LFT)
$
f(M,K)=M_{11}+M_{12}K(I-M_{22}K)^{-1}M_{21}$ is given by the state-space formula
\begin{align} \label{eq:LFT}
f(M,K)=\ll \ba{cc|c}
A+B_2D_KC_2 & B_2C_K & B_1+B_2D_KD_{21} \\
B_KC_2 & A_K& B_KD_{21} \\ \hline
C_1+D_{12}D_KC_2 & D_{12}C_K & D_{11}+D_{12}D_{K}D_{21}
\ea \rr.
\end{align}

\end{lemma}
\begin{proof}
The proof is standard, see for example \cite[pp. 179]{Zhou} and the references therein.
\end{proof}

\begin{proof}[Proof of Theorem \ref{theorem:3}] 
Consider again the optimal control problem \eqref{eq:chap4_opt_controll}:
\begin{equation}\label{eq:A1} 
\ba{rl}
\underset{K}{\text{minimize}} & \|P_{11}+P_{12}K(I-P_{22}K)^{-1}P_{21} \|^{2} \\
\text{subject to \ \ } & K\in \mc{I}(\mc{P}) \\
& K \text{ stabilizing}.
\ea
\end{equation}
Let $v_1^{*}$ be the optimal value of \eqref{eq:A1}. Consider, on the other hand the optimization problem:
\begin{equation} \label{eq:A3}
\ba{rl}
\underset{Q}{\text{minimize}} & \|P_{11}+P_{12}Q \|^{2} \\
\text{subject to \ \ } & Q \in \mc{I}(\mc{P}). 
\ea
\end{equation}
Let $v_2^{*}$ be the optimal value of \eqref{eq:A3}. Recall that the optimal solution $Q^{*}$ of \eqref{eq:A3} was obtained in Lemma \ref{lemma:A4} as \eqref{eq:A2}.
We note that if $K^{*}$ is an optimal solution to \eqref{eq:A1} then the corresponding $\bar{Q}:=K^{*}(I-P_{22}K^{*})^{-1}P_{21}$ is feasible for \eqref{eq:A3}. Hence $v_2^{*} \leq v_1^{*}$. We will show that the controller in \eqref{eq:thm2} is optimal by showing that $\bar{Q}=Q^{*}$ (so that $v_1^{*}=v_2^{*}$). We will also show that $K^{*} \in \mc{I}(\mc{P})$ and is internally stabilizing. Since it achieves the lower bound $v_2^{*}$ and is internally stabilizing, it \emph{must} be optimal.

Given $K^{*}$, one can evaluate $\bar{Q}:=K^{*}(I-P_{22}K^{*})^{-1}P_{21}$. To do so we use $K^{*}$ as per \eqref{eq:thm2} and 
$$
M=\ll \ba{cc}
0 & I \\  P_{21} & P_{22}
\ea \rr =
\ll \ba{c|cc}
A & F & B \\ \hline
0 & 0 & I \\
I & 0 & 0
\ea \rr
$$
and use the interconnection formula \eqref{eq:LFT}
to obtain:
$$
\bar{Q}=\ll \ba{cc|c}
A+BC_Q\Pi_1 & BC_{Q}(\Pi_2-\Pi_1C_\Phi) & F \\
B_{\Phi} &  A_{\Phi}-B_{\Phi}C_{\Phi} & 0 \\ \hline
C_{Q}\Pi_1 & C_{Q}(\Pi_2-\Pi_1C_{\Phi}) & 0
\ea \rr.
$$

Recall that $Q^{*}$ given by \eqref{eq:A2} is the optimal solution to \eqref{eq:6a} (which constitutes a lower bound to the problem we are trying to solve). We are trying to show that it is achievable by explicitly producing $K^{*}$ such that $\bar{Q}:=K^{*}(I-P_{22}K^{*})^{-1}P_{21}$ and $\bar{Q}=Q^{*}$, thereby proving the optimality of $K^{*}$.
 
While $Q^{*}$ in \eqref{eq:A2} and $\bar{Q}$ obtained above appear different at first glance, their state-space realizations are actually equivalent modulo a coordinate transformation. Recall that $\Pi_2$ is a matrix (composed of standard unit vectors) that spans the orthogonal complement of the column span of $\Pi_1$. As a result the matrix $\ll \ba{cc} \Pi_1 & \Pi_2 \ea \rr$ is a permutation matrix. Define the matrices
\begin{align*}
\Lambda:=\ll \ba{cc} \Pi_1 & \Pi_2  \ea \rr \ll \ba{cc}
I & -C_{\Phi} \\
0 & I
\ea \rr \text{ , \ } &
\Lambda^{-1}= \ll \ba{cc}
I & C_{\Phi} \\
0 & I
\ea \rr \ll \ba{c} \Pi_1^{T} \\ \Pi_2^{T}  \ea \rr.
\end{align*}

Note that $\Lambda$ is a square, invertible matrix.
 Changing state coordinates on $Q^{*}$ using $\Lambda$ via:
\begin{align*}
\mathbf{A} & \mapsto \Lambda^{-1} \mathbf{A}\Lambda \\
\Pi_1F & \mapsto \Lambda^{-1} \Pi_1F \\
C_Q & \mapsto C_{Q} \Lambda
\end{align*}
along with the relations $R\Pi_2\Pi_2^{T}+\Pi_1^{T}=R$, $R\mathbf{A}-BC_Q=AR$, and $AR\Pi_1=A$, we see that the transformed realization of $Q^{*}$ is equal to the realization of $\bar{Q}$, and hence $Q^{*}=\bar{Q}$.

Using \eqref{eq:1} for the open loop, \eqref{eq:thm2} for the controller and the LFT formula \eqref{eq:LFT} to compute the closed loop map, one obtains that the closed-loop state transition matrix is given by
$$
\ll \ba{cc}A+BC_Q\Pi_1 & BC_{Q}(\Pi_2-\Pi_1C_\Phi) \\
B_{\Phi} &  A_{\Phi}-B_{\Phi}C_{\Phi} \ea \rr \cong \mathbf{A}.
$$ 
By Lemma \ref{lemma:diagstab}, the closed loop is internally stable.

By the column concatenation formula and \eqref{eq:lem_inc} we have
\begin{align*} 
\ll \ba{c}\Phi \\ K_\Phi\ea \rr &= \ll \ba{c|c}
A_{\Phi} & B_{\Phi} \\ \hline
C_\Phi & I \\
C_Q\Pi_2 & C_Q \Pi_1
\ea \rr 
\end{align*}
Using the state-space coprime factorization formula \cite[pp. 52]{Zhong} it is straightforward to verify that for the expression for $K^{*}$ in \eqref{eq:thm2},
$K^{*}=K_{\Phi}\Phi^{-1}$. Since by Lemma \ref{lemma:inc_alg} both $\Phi \in \mc{I}(\mc{P}) $ and $K_{\Phi} \in \mc{I}(\mc{P})$, we have $K^{*} \in \mc{I}(\mc{P})$. 
\end{proof}

\begin{proof}[Proof of Theorem \ref{theorem:4}] 
By the column concatenation formula and \eqref{eq:lem_inc} we have
\begin{align*} 
\ll \ba{c}\Phi \\ K_\Phi\ea \rr &= \ll \ba{c|c}
A_{\Phi} & B_{\Phi} \\ \hline
C_\Phi & I \\
C_Q\Pi_2 & C_Q \Pi_1
\ea \rr 
\end{align*}
Using the state-space coprime factorization formula \cite[pp. 52]{Zhong} to evaluate $K_\Phi \Phi^{-1}$ we see that $K^{*}=K_\Phi \Phi^{-1}$. This directly gives the first expression in the statement of the theorem. 
The second expression is a simple manipulation of the first expression.
\end{proof}

\section{Conclusions}
In this paper we provided a state-space solution to the problem of computing an $\mc{H}_{2}$-optimal decentralized controller for a poset-causal system. We introduced a new decomposition technique that enables one to separate the decentralized problem into a set of centralized problems. We gave explicit state-space formulae for the optimal controller and provided degree bounds on the controller. We illustrated our technique with a numerical example. Our approach also enabled us to provide insight into the structure of the optimal controller. We introduced a pair of transfer functions $(\Phi, \Gamma)$ and showed that they were intimately related to the prediction of the state.
In future work, it would be interesting to attempt to apply this decomposition technique to a wider class of decentralization structures.

\bibliography{h2_journal_SUBMISSION}

\appendix
\subsection{Numerical Example} \label{sec:num}
In this section, we consider a numerical example for the poset shown in Fig. \ref{fig:2.1}(d). The system has one state and one input per subsystem, and we synthesize the optimal controller. The data for the is the same as in Example \ref{example:2} with the matrices $A,B,C,D$ as given in \eqref{eq:comp_example}.
%
For this problem the relevant matrices that are used in constructing the controller are:
$$ 
\ba{rl}
\Pi_1&=\ll \ba{cccc} E_{1} & E_{5} & E_{7} &E_{9} \ea \rr \\ \\
R&=\ll \ba{cccc|cc|cc|c}         
1 & 0 & 0& 0& 0 & 0 & 0 & 0 & 0  \\
 0& 1& 0& 0&1 &0 &0 &0 &0 \\
 0& 0&1& 0& 0 &0 &1 &0 &0\\
0& 0& 0& 1& 0 &1 &0 &1 &1
 \ea
 \rr \\ \\ 

\Pi_2 &=\ll \ba{ccccc} E_{2} & E_{3} & E_{4} & E_{6} &E_8 \ea \rr.
\ea
$$ 
(Recall that $E_i$ is the $9 \times 1$ $i^{th}$ unit vector.)
Note that $\down 1 = \l 1,2,3,4 \r, \down 2 = \l 2,4 \r, \down 3 = \l 3,4 \r, \down 4 = \l 4 \r$. Accordingly, the Riccati subproblems that we need to solve are given by $\text{Ric}(\down j)$ for $j=1,2,3,4$. 
Upon solving these, we obtain
\begin{align*}
K(\down 1, \down 1)&= \ll \ba{cccc}  .7175 & .3515 & .3616 & -.0751\\
-.9671 & .9575 & .1827 &.1033 \\
-1.0306 & .2045 & 1.0312 & .0814 \\
.6337 & -.7902 & -.8121 & .8935  \ea \rr & K(\down 2, \down 2)&=\ll \ba{cc} 
         1.0237 & .0990 \\
-.8011 & .9001
 \ea \rr \\
K(\down 3, \down 3)&=\ll \ba{cc} 1.0960 & .0792 \\ -.8226 & .9019  \ea \rr& K(\down 4, \down 4)&=.9050
.
\end{align*}
From these, it is possible to construct $A(\down j, \down j)-B(\down j, \down j)K(\down j, \down j)$ for $j \in \left\{ 1, 2, 3,4 \right\}$, and from that construct $\mathbf{A}=\text{diag}(A(\down j, \down j)-B(\down j, \down j)K(\down j, \down j))$ given by
$$
\mathbf{A}=\ll \ba{ccccccccc}
-1.2175 & -.3515 & -.3616 & .0751&0&0&0&0&0\\
-.7505 & -1.5589 & -.5443 & -.0282&0&0&0&0&0\\
-.6869&-.5560&-1.5927&-.0064&0&0&0&0&0\\
-.3535 & -1.7232 & -1.7634 & -1.1032&0&0&0&0&0\\
0&0&0&0&-1.2737 & -.0990&0&0&0\\
0&0&0&0&-1.2226 & -1.0991&0&0&0\\
0&0&0&0&0&0&-1.2969 & -.0792&0\\
0&0&0&0&0&0&-1.2733 & -1.0811&0\\
0&0&0&0&0&0&0&0&-1.0050\\
\ea \rr.
$$
Using \eqref{eq:matrix_defs} one readily obtains $A_{\Phi}$ and $B_{\Phi}$ to be
\begin{align*}
A_{\Phi}&= \ll \ba{ccccc} 
-1.5589 & -.5443 & -.0282 & 0 & 0 \\
-.5560 & -1.5927 & -.0064 & 0 & 0 \\
-1.7232 & -1.7634 & -1.1032 & 0 & 0 \\
0 & 0 & 0 & -1.0991 & 0 \\
0 & 0 & 0 & 0 & -1.0811
\ea \rr &
B_{\Phi}&= \ll \ba{cccc}
-.7505 & 0 & 0 & 0  \\
-.6869 & 0 & 0 & 0  \\
-.3535 & 0 & 0& 0\\
0 &-1.2226 &0&0 \\
0 & 0 & -1.2733 & 0
\ea \rr.
\end{align*}
Note that given $K(\down 2, \down 2)$ ( a $2 \times 2$ matrix) one needs to construct $\hat{K}(\down 2, \down 2)$ (a $4 \times 4$ matrix) by zero padding. For instance, we have
$$
\hat{K}(\down 2, \down 2)=\ll \ba{cccc}
0         &    0        &     0     &       0 \\
             0       &  1.0237  &           0    &     .0990 \\
             0     &        0   &          0    &         0 \\
             0     &    -.8011    &         0    &     .9001
 \ea \rr.
$$
From these,  one constructs $C_Q$ using \eqref{eq:matrix_defs} to be
$$
C_Q=\ll \ba{ccccccccc}
.7175 & .3515 & .3616 & -.0751&0&0&0&0&0 \\
-.9671 & .9575 & .1827 &.1033& 1.0237&.0990&0&0&0 \\
-1.0306 & .2045 & 1.0312 & .0814&0&0&1.0960&.0792 &0 \\
.6337 & -.7902 & -.8121 & .8935&-.8011&..9001&-.8226&.9019&.9050 
\ea \rr.
$$
We use these quantities to obtain the controller $K^{*}$ using formula \eqref{eq:thm2}. 
The controller 
$$
K=\ll \ba{c|c} A_K & B_K \\ \hline C_K & D_K \ea \rr,
$$
has the following realization: 
\begin{align*}
A_K&= \ll \ba{ccccc}   -1.5589     &        -.5443       & -.0292 & 0&     0 \\
         -.5560 & -1.5927 & -.0064      &       0      &       0 \\
        -1.7232 & -1.7634 & -1.1032     &        0    &         0 \\
          1.2226   &          0     &        0   &      -1.0991   &          0 \\
             0     &     1.2733     &        0    &         0     &    -1.0811 \ea \rr  
              \end{align*}
 \begin{align*}    
B_K&= \ll \ba{cccc}
-.7505     &        0      &       0     &        0 \\
         -.6869   &          0   &          0     &        0 \\
          -.3535   &          0   &          0     &        0 \\
             0     &    -1.2226    &         0    &         0 \\
             0      &       0   &      -1.2733     &        0
 \ea \rr  
 \end{align*}
 \begin{align*}     
C_K&=   \ll \ba{ccccc}
          -.7175     &     -.3616 &       .0751     &        0    &         0 \\
         .0662      &    -.1827   &       -.1033    &     -.0991   &          0 \\
         -2.045   & .0648 & -.0814 &         0       &   -.0792 \\
         -.0109 & -.0106 & .0115 & .0489 & .0031
 \ea \rr 
  \end{align*}
 \begin{align*}    
D_K&= \ll \ba{cccc} -.7175      &       0      &       0     &        0 \\
          .9671   &       -1.0237   &         0     &        0\\
          1.0306   &          0      &   -1.0960    &         0\\
         -.6337  &        0.8011    &      0.8226   &       -0.9050 \ea \rr. 
\end{align*} 
Note that the optimal controller is of degree $5$. This matches the bound obtained in Corollary~\ref{cor:degbounds} exactly. Note also that the matrix $D_K$ is in the incidence algebra (and so is the controller $K$ itself, as can be verified from the transfer function). Finally, this controller can be verified to be stabilizing. Let $h_{\text{open}}, h_{\text{centralized}}, h_{\text{decentralized}}$ be the open loop, optimal centralized closed loop and optimal decentralized closed loop $\mc{H}_2$ norms. We obtain the following values:
\begin{align*}
h_{\text{open}}&=31.6319 \\
h_{\text{centralized}}&=2.3197 \\
h_{\text{decentralized}}&=2.8280.
\end{align*}

\end{document}